\documentclass[11pt]{article}

\usepackage[cmtip,arrow]{xy}
\usepackage{amsmath,amssymb,enumerate,pb-diagram,pb-xy}
\def\today{\number\day .\number\month .\number\year}
\usepackage[applemac]{inputenc}

\def\1{{\bf 1}}
\def \a{{\mathfrak a}}
\def \Ad{{\rm Ad}}
\def \al{\alpha}

\def \b{{\mathfrak b}}
\def \bs{\backslash}
\def \C{{\mathbb C}}

\def \CE{{\cal E}}

\def \CN{{\cal N}}
\def \CO{{\cal O}}

\def \CV{{\cal V}}
\def\det{\operatorname{det}}

\def\e{\emph}
\def \End{{\rm End}}
\def\Ext{{\rm Ext}}
\def \g{{\mathfrak g}}
\def \ga{\gamma}
\def \Ga{\Gamma}
\def \h{{\mathfrak h}}
\def \Hom{{\rm Hom}}

\def \k{{\mathfrak k}}
\def \la{\lambda}
\def\l{{\mathfrak l}}

\def \m{{\mathfrak m}}

\def\mathqed{\vspace{-25pt}\newline \vspace{0pt}\qed\\ \vspace{25pt}}
\def \n{{\mathfrak n}}
\def \N{{\mathbb N}}
\def \ord{{\rm ord}}
\def \p{{\mathfrak p}}
\def \ph{\varphi}

\def \prf{{\bf Proof: }}

\def\qed{\hfill $\square$}
\def \R{{\mathbb R}}
\def \Re{\operatorname{Re}}

\def \t{{\mathfrak t}}
\def \tr{{\rm tr}}

\def \vol{{\rm vol}}

\def \Z{{\mathbb Z}}
\def \({\left(}
\def \){\right)}
\def \={{\ =\ }}

\newtheorem{theorem}{Theorem}[section]

\newtheorem{lemma}[theorem]{Lemma}

\newtheorem{proposition}[theorem]{Proposition}

\begin{document}

\pagestyle{myheadings} \markright{HOLOMORPHIC TORSION}

\title{Holomorphic torsion and closed geodesics}
\author{Anton Deitmar\\ \ \\
J. Fixed Point Theory Appl. 9, 25-53 (2011)}
\date{}
\maketitle

$$ $$

\tableofcontents

\newpage

\section*{Introduction}
In his seminal paper \cite{Fried86}, David Fried expressed the analytic torsion of a hyperbolic manifold as a special value of a zeta function.
Later in \cite{Fried88} he did the same for the holomorphic torsion of a complex hyperbolic manifold.
In the paper \cite{MS-tors}, based on \cite{MS-eta},  Moscovici and Stanton generalized Fried's work on analytic torsion to locally symmetric spaces of fundamental rank one. 

In the current paper we extend Fried's result on holomorphic torsion to arbitrary hermitian locally symmetric spaces.
In the papers \cite{Fried86,Fried88}, the torsion numbers are described in the form $R(0)c$, where $R$ is a geometrically defined zeta function and $c$ only depends on the volume and the universal cover of the space.
In this paper we interpret $c$ as the reciprocal of the $L^2$-torsion.
The main result of this paper can also be viewed as a geometric analogue of Lichtenbaum's conjecture \cite{Licht}.
See also \cite{Den}.

Fried's first paper on analytic torsion \cite{Fried86}, which only treated hyperbolic spaces, has been generalized to spaces of fundamental rank one by Moscovici and Stanton in \cite{MS-tors} and by the author to all spaces of positive fundamental rank \cite{D-Hitors}.
It then became clear, that the analogous treatment of holomorphic torsion requires a different approach.
The basic problem is to find suitable test functions for the trace formula to read off the analytic continuation of Selberg-type zeta functions.
To explain this, we fix a semisimple Lie group $G$ together with a cocompact lattice $\Ga\subset G$.
By integration of the right regular representation, a function $f\in C_c^\infty(G)$ induces an operator $R(F)$
on the Hilbert space $L^2(\Ga\bs G)$.
Having a smooth kernel, the operator $R(f)$ is trace class and its trace equals the integral over the diagonal of the kernel.
Computing this integral, one gets an identity, known as the trace formula, expressing the (spectral) trace $\tr R(f)$ in terms of a (geometric) sum of orbital integrals.
For suitable test functions $f$, the geometric side gives Selberg-type zeta functions and the comparison with the spectral side yields analytic continuation of the zeta function.
To find such test functions, one can either look for test functions one knows the spectral side of, or for ones with given geometric trace.
The first, ``spectral'' approach is the classical one.
Usually one takes functions which are given by the functional calculus of a nice invariant operator like a Laplace operator.
Fried and Moscovici/Stanton used heat kernels to this end.
For higher rank groups, however, this approach cannot alway distinguish geometric contributions from different Cartan subgroups, which is necessary to achieve analytic continuation.
At this point, Juhl's habilitation thesis \cite{Ju} gave a new idea.
It is the first treatment systematically to employ the second, geometric approach.
It is detailed in the book \cite{Juhl}.
This technique only works in rank one situations, but it inspired the author to \e{combine} the two approaches in the following way.
First, the geometric technique of Juhl is used to reduce the rank by one, and then a spectral construction is used to give a kernel on a Levi-subgroup which eliminates all unwanted contributions.
The result, which is given in \cite{lefschetz}, is a Lefschetz formula isolating the geometric contributions of a single Cartan subgroup.

This result is more general then what we need in this paper, as it is valid for all semisimple Lie groups.
We will specialize it here to groups attached to hermitian symmetric spaces which requires a detailed knowledge of the latter.

\section{Notation}
This paper depends heavily on \cite{lefschetz}.
We take over all notation from that paper.
For instance, $G$ will denote a connected semisimple Lie group with finite center, and $K$ will be a maximal compact subgroup.
The attached locally symmetric space is denoted by $X=G/K$.
For every discrete, torsion-free subgroup $\Ga\subset G$ the space $\Ga\bs X$ is a locally symmetric space.
We will only be interested in the compact case, so $\Ga\bs X$, or equivalently, $\Ga\bs G$ will be compact.
In this case $\Ga$ is said to be a \emph{uniform lattice} in $G$.

Let $L$ be a Lie group. An element $x$ of
$L$ is called \emph{neat}\index{neat} if  for every
 finite dimensional representation $\eta$ of $L$ the linear map $\eta(x)$ 
 has no nontrivial root of unity as an eigenvalue. 
A subset $A$ of $L$ is called neat if each of its members are. Every neat subgroup is torsion free modulo the center of $L$.
Every arithmetic group has a subgroup of finite index which is neat
\cite{Bor}.

\begin{lemma}
Let $x\in L$ be semisimple and neat. Let $L_x$ denote the centralizer of $x$
in
$L$. Then for each $k\in\N$ the connected components of $L_x$ and $L_{x^k}$ coincide.

\end{lemma}

\prf
It suffices to show that the Lie algebras coincide. The Lie algebra of $L_x$ is just
the fixed space of $\Ad(x)$ in Lie$(L)$. Since $\Ad(x)$ is semisimple and does not have a
root of unity for an eigenvalue this fixed space coincides with the fixed space of
$\Ad(x^k)$. The claim follows.
\qed

Let $\Ga\subset G$ denote a cocompact discrete subgroup which is
 neat. Since $\Ga$ is torsion free it acts fixed point free
on the contractible space $X$ and hence $\Ga$ is the fundamental
group of the Riemannian manifold
$$
X_{\Ga} = \Ga \bs X = \Ga \bs G/K
$$\index{$X_\Ga$}
 it follows that we have a
canonical bijection of the homotopy classes of loops:
$$
[S^1 : X_{\Ga} ] \to \Ga / {\rm conjugacy}.
$$
For a given class $[\ga]$ let $X_\ga$\index{$X_\ga$} denote the
union of all closed geodesics in the corresponding class in $[S^1
: X_\Ga ]$. Then $X_\ga$ is a smooth submanifold of $X_{\Ga_H}$
\cite{DKV}, indeed, it follows that
$$
X_\ga\ \cong\ \Ga_\ga \bs G_\ga/K_\ga,
$$
where $G_\ga$ and $\Ga_\ga$ are the centralizers of $\ga$ in $G$
and $\Ga$ and $K_\ga$ is a maximal compact subgroup of $G_\ga$.
Further all closed geodesics in the class $[\ga]$ have the same
length $l_\ga$.\index{$l_\ga$}

\begin{lemma}\label{neatXgamma}
For $\ga\in\Ga$ and $n\in\N$ we have $X_{\ga^n}=X_\ga$.
\end{lemma}

\prf By the last lemma we have for the connected components, $G_\ga^0=G_{\ga^n}^0$.
By definition one has that $X_{\ga^n}$ is a subset of $X_\ga$. Since both are connected
submanifolds of $X_\Ga$ they are equal if their dimensions are the same.
We have
\begin{eqnarray*}
\dim\, X_\ga &=& \dim\, G_\ga/K_\ga\\
&=& \dim\, G_\ga^0/K_\ga^0\\
&=& \dim\, G_{\ga^n}^0/K_{\ga^n}^0\\
&=& \dim\, G_{\ga^n}/K_{\ga^n}\\
&=&\dim\, X_{\ga^n}. 
\end{eqnarray*}
\mathqed

\section{Meromorphic continuation}
Recall that $G$ is a connected semisimple Lie
group with finite center. Fix a maximal compact subgroup
$K$ with Cartan involution $\theta$.

\subsection{The zeta function} \label{constzeta}

Fix a $\theta$-stable Cartan subgroup $H$ of split rank 1. Note
that such a $H$ doesn't always exist. It exists only if the
absolute ranks of $G$ and $K$ satisfy the relation:
$$
{\rm rank}\ G -{\rm rank}\ K \leq 1.
$$
This certainly holds if $G$ has a compact Cartan subgroup or if
the real rank of $G$ is one. In the case ${\rm rank}\ G ={\rm
rank}\ K$, i.e., if $G$ has a compact Cartan there will in general
be several $G$-conjugacy classes of split rank one Cartan
subgroups. In the case ${\rm rank}\ G -{\rm rank}\ K = 1$,
however, there is only one. The number $FR(G):={\rm rank}\ G
-{\rm rank}\ K $ is called the {\it fundamental
rank}\index{fundamental rank} of $G$.

Write $H=AB$ where $A$ is the connected split component and
$B\subset K$ is compact. Choose a parabolic subgroup $P$
with Langlands decomposition $P=MAN$. Then $K_M =K\cap M$
is a maximal compact subgroup of $M$. Let
$\bar{P}=MA\bar{N}$ be the opposite parabolic. As in the
last section let $A^-=\exp(\a_0^-)\subset A$ be the
negative Weyl chamber of all $a\in A$ which act
contractingly on $\n$. Fix a finite dimensional
representation $(\tau ,V_\tau)$ of $K_M$.

Let $H_1\in \a_0^-$ be the unique element with $B(H_1)=1$.

Let $\Ga\subset G$ be a cocompact discrete subgroup which is
 neat. An element $\ga \neq 1$ will be called {\it
primitive}\index{primitive element} if $\tau \in \Ga$ and $\tau^n
=\ga$ with $n\in \N$ implies $n=1$. Every $\ga \neq 1$ is a power
of a unique primitive element. Obviously primitivity is a property
of conjugacy classes. Let $\CE_P^p(\Ga)$ denote the subset of
$\CE_P(\Ga)$ consisting of all primitive classes. Recall the
length $l_\ga$ of any geodesic in the class $[\ga]$. If $\ga$ is
conjugate to $am\in A^-M_{ell}$ then $l_\ga=l_a$, where $l_a=|\log
a|=\sqrt{B(\log a)}$. Let $A_\ga$ be the connected split component
of $G_\ga$, then $A_\ga$ is conjugate to $A$.

We say that $M$ is \emph{orientation preserving}, if $M$ acts by orientation preserving maps on the manifold $M/K_M$, where $K_M=K\cap M$.
For a complex vector space $V$, on which $A$ acts linearly, and $\la\in \a^*$, we write $V^\la$ for the generalized $\la$-Eigenspace, i.e.,
$$
V^\la\=\{v\in V: (a-\la)^kv=0\text{ for some }k\in\N\}.
$$

In this section we are going to prove the following

\begin{theorem}\label{genSelberg} Suppose that $\tau$ is in
the image of the restriction map $res_{K_M}^M$ or that $M$
is orientation preserving. Let $\Ga$ be  neat and
$(\omega ,V_\omega)$ a finite dimensional unitary
representation of $\Ga$. For $\Re(s)>>0$ define the {\it
generalized Selberg zeta function}:\index{Selberg zeta
function}
$$
Z_{P,\tau,\omega}(s) \= \prod_{[\ga]\in {\cal E}_P^p(\Ga)}
\prod_{N\geq 0} \det\left(1-e^{-sl_\ga}\ga \left|
\begin{array}{c}V_\omega \otimes V_\tau \otimes
S^N(\n)\end{array}\right.\right)^{\chi(A_\ga\bs X_\ga)},
$$
where $S^N(\n)$ denotes the $N$-th symmetric power of the space
$\n$ and $\ga$ acts on $V_\omega \otimes V_\tau \otimes S^N({\n})$
via $\omega(\ga) \otimes \tau(b_\ga) \otimes Ad^N(a_\ga b_\ga)$,
here $\ga \in \Ga$ is conjugate to $a_\ga b_\ga \in A^-B$.
Finally, $\chi(A_\ga\bs X_\ga)$ is the Euler-characteristic of $A_\ga\bs X_\ga$.

Then $Z_{P,\tau,\omega}$ has a meromorphic continuation to the
entire plane. The vanishing order of $Z_{P,\tau ,\omega}(s)$ at a
point $s=\la (H_1)$, $\la \in \a^*$, is
{\small
$$
m(s)=(-1)^{\dim N}\sum_{\pi \in \hat{G}}N_{\Ga ,\omega}(\pi)
\sum_{p,q}(-1)^{p+q} \dim
\left(\begin{array}{c}H^q(\n,\pi_K)^\la\otimes \wedge^p\p_M \otimes
V_{\breve{\tau}}\end{array}\right)^{K_M}.
$$}
Further, all poles and zeroes of the function $Z_{P,\tau
,\omega}(s+|\rho_0|)$ lie in $\R \cup i\R$.
\end{theorem}

\prf
We apply the Lefschetz Theorem (Thm 4.2.1) of \cite{lefschetz} with the test function $\ph(a)=l_a^{j+1}\,e^{-sl_a}$ for $j\in\N$ and $s\in\C$ with $j,\Re(s)$ large enough.
The geometric side of the Lefschetz formula is $\(\frac{\partial}{\partial s}\)^j\frac{Z'}{Z}(s_)$, and by the Lefschetz formula this equals $\(\frac\partial{\partial s}\)^j\sum_{s_0\in\C}\frac{m(s_0)}{s-s_0}$.
The theorem follows.
\qed

\subsection{The functional equation}

If the Weyl group $W(G,A)$ is nontrivial, then is has
order two. To give the reader a feeling of this condition consider
the case $G=SL_3(\R)$. In that case the Weyl group $W(G,A)$ is
trivial. On the other hand, consider the case when the fundamental
rank of $G$ is $0$; this is the most interesting case to us since
only then we have several conjugacy classes of splitrank-one
Cartan subgroups. 
Here the \e{splitrank} of a Cartan subgroup $H=AB$ is the dimension of $A$, where $B$ is compact and $A$ is a split torus.
In that case it follows that the dimension of
all irreducible factors of the symmetric space $X=G/K$ is even,
hence the point-reflection at the point $eK$ is in the connected
component of the group of isometries of $X$. This reflection can
be thought of as an element of $K$ which induces a nontrivial
element of the Weyl group $W(G,A)$. So we see that in this
important case we have $|W(G,A)|=2$.

Let $w$ be the nontrivial element. It has a representative in $K$
which we also denote by $w$. Then $wK_Mw^{-1}=K_M$ and we let
$\tau^w$ be the representation given by
$\tau^w(k)=\tau(wkw^{-1})$. It is clear from the definitions that
$$
Z_{P,\tau,\omega}\= Z_{\bar{P},\tau^w,\omega}.
$$
We will show a functional equation for $Z_{P,\tau,\omega}$.
This
needs some preparation. Assume $G$ admits a compact Cartan
$T\subset K$, then a representation $\pi\in\hat{G}$ is called
\emph{elliptic} if $\Theta_\pi$ is nonzero on the compact Cartan. Let
$\hat{G}_{ell}$ be the set of elliptic elements in $\hat{G}$ and
denote by $\hat{G}_{ds}$ the subset of discrete series
representations. Further let $\hat{G}_{lds}$ denote the set of all
discrete series and all limits of discrete series representations.
In Theorem \ref{genSelberg} we have shown that the vanishing order
of $Z_{P,\tau,\omega}(s)$ at the point $s=\mu (H_1)$, $\mu\in\a^*$
is
$$
(-1)^{\dim N}\sum_{\pi\in\hat{G}} N_{\Ga ,\omega}(\pi) m(\pi ,\tau
,\mu),
$$
where
$$
m(\pi ,\tau ,\mu) \= \sum_{p,q} (-1)^{p+q} \dim
\left(\begin{array}{c}H^q(\n ,\pi_K)(\mu)\otimes \wedge^p\p_M\otimes
V_{\breve{\tau}}\end{array}\right)^{K_M}.
$$

A \e{standard representation} is a representation parabolically induced from a tempered representation, see \cite{Knapp}, p.383. 
Any character $\Theta_\pi$ for $\pi\in\hat{G}$ is an integer
linear combination of characters of standard representations (loc. cit.). From
this it follows that for $\pi\in\hat{G}$ the character restricted
to the compact Cartan $T$ is
$$
\Theta_\pi \mid_T \= \sum_{\pi'\in\hat{G}_{lds}} k_{\pi ,\pi'}
\Theta_{\pi'} \mid_T,
$$
with integer coefficients $k_{\pi ,\pi'}$.

\begin{lemma}
There is a $C>0$ such that for $\Re (\mu(H_1))<-C$ the order of
$Z_{P,\tau,\omega}(s)$ at $s=\mu(H_1)$ is
$$
(-1)^{\dim N}\sum_{\pi\in\hat{G}_{ell}} N_{\Ga ,\omega}(\pi)
\sum_{\pi'\in\hat{G}_{lds}} k_{\pi ,\pi'} m(\pi' ,\tau ,\mu).
$$
\end{lemma}

\prf For any $\pi\in\hat{G}$ we know that if $\Theta_\pi
|_{AB}\neq 0$ then in the representation of $\Theta_\pi$ as
linear combination of standard characters there must occur
lds-characters and characters of representations $\pi_{\xi ,\nu}$
induced from $P$. Since $\Theta_{\pi_{\xi
,\nu}}=\Theta_{\pi_{^w\xi ,-\nu}}$ any contribution of $\pi_{\xi
,\nu}$ for $\Re (s) <<0$ would also give a pole or zero of
$Z_{P,\tau^w,\omega}$ for $\Re(s)>>0$. In the latter region we do
have an Euler product, hence there are no poles or zeroes. \qed

Now consider the case $FR(G)=1$, so there is no compact Cartan,
hence no discrete series.

\begin{theorem} Assume that the fundamental rank of $G$ is $1$,
then there is a polynomial $P$ of degree $\leq \dim G+\dim N$ such
that
$$
Z_{P,\tau,\omega}(s) \= e^{P(s)} Z_{P,\tau^w,\omega}(2| \rho_0|
-s).
$$
\end{theorem}

\prf By Theorem \ref{genSelberg} the functions $Z_{P,\tau,\omega}(s)$ and $Z_{P,\tau^w,\omega}(2| \rho_0|
-s)$ have the same zeros and poles.
Further, they are both of finite genus, hence the claim.
\qed

Now assume $FR(G)=0$ so there is a compact Cartan subgroup $T$. As
Haar measure on $G$ we take the Euler-Poincar\'e measure. The sum
in the lemma can be rearranged to
$$
(-1)^{\dim N} \sum_{\pi'\in\hat{G}_{lds}}  m(\pi' ,\tau
,\mu)\sum_{\pi\in\hat{G}_{ell}} N_{\Ga ,\omega}(\pi) k_{\pi
,\pi'}.
$$
We want to show that the summands with $\pi'$ in the limit of the
discrete series add up to zero. For this suppose $\pi'$ and
$\pi''$ are distinct and belong to the limit of the discrete
series. Assume further that their Harish-Chandra parameters agree.
By the Paley-Wiener theorem \cite{CloDel} there is a smooth
compactly supported function $f_{\pi' ,\pi''}$ such that for any
tempered $\pi\in\hat{G}$:
$$
\tr \pi(f_{\pi' ,\pi''}) \= \left\{ \begin{array}{cl} 1& {\rm if}\
\pi =\pi'\\
                        -1  &{\rm if}\ \pi =\pi''\\
                        0   & {\rm else.}
                    \end{array}\right.
$$
Plugging $f_{\pi' ,\pi''}$ into the trace formula one gets
$$
\sum_{\pi\in\hat{G}_{ell}} N_{\Ga ,\omega}(\pi)k_{\pi ,\pi'}
\=\sum_{\pi\in\hat{G}_{ell}} N_{\Ga ,\omega}(\pi)k_{\pi ,\pi''},
$$
so that in the above sum the summands to $\pi'$ and $\pi''$ occur
with the same coefficient. Let $\pi_0$ be the induced
representation whose character is the sum of the characters of the
$\pi''$, where $\pi''$ varies over all lds-representations with
the same Harish-Chandra parameter as $\pi'$. Then for $\Re
(\mu(H_1))<-C$ we have $m(\pi_0,\tau,\mu)=0$. Thus it follows that
the contribution of the limit series vanishes.

Plugging the pseudo-coefficients \cite{Lab} of the discrete series
representations into the trace formula gives for
$\pi\in\hat{G}_{ds}$:
$$
\sum_{\pi'\in\hat{G}_{ell}} k_{\pi' ,\pi} N_{\Ga ,\omega}(\pi') \=
\dim \omega (-1)^{\frac{\dim X}{2}}\chi(X_\Ga) d_{\pi},
$$
where $d_{\pi}$ is the formal degree of $\pi$.

The infinitesimal character $\la$ of $\pi$ can be viewed as an
element of the coset space $(\t^*)^{reg}/W_K$. So let $J$ denote
the finite set of connected components of $(\t^*)^{reg}/W_K$, then
we get a decomposition $\hat{G}_{ds} = \coprod_{j\in J}
\hat{G}_{ds,j}$.

In the proof of Theorem \ref{genSelberg} we used the Hecht-Schmid character
formula to deal with the global characters. On the other hand it
is known that global characters are given on the regular set by
sums of toric characters over the Weyl denominator. So on $H=AB$
the character $\Theta_\pi$ for $\pi\in\hat{G}$ is of the form $\CN
/D$, where $D$ is the Weyl denominator and the numerator $\CN$ is
of the form
$$
\CN (h) \= \sum_{w\in W(\t,\g)} c_w h^{w\la},
$$
where $\la\in\h^*$ is the infinitesimal character of $\pi$.
Accordingly, the expression $m(\pi ,\tau ,\mu)$ expands as a sum
$$
m(\pi ,\tau ,\mu) \= \sum_{w\in W(\t ,\g)} m_w(\pi ,\tau ,\mu).
$$

\begin{lemma}
Let $\pi ,\pi'\in \hat{G}_{ds,j}$ with infinitesimal characters
$\la ,\la'$ which we now also view as elements of $(\h^*)^+$, then
$$
m_w(\pi ,\tau ,\mu) \= m_w(\pi',\tau ,\mu +w(\la'-\la)|_\a).
$$
\end{lemma}

\prf In light of the preceding it suffices to show the following:
Let $\tau_\la, \tau_{\la'}$ denote the numerators of the global
characters of $\pi$ and $\pi'$ on $\h^+$. Write
$$
\tau_\la (h) \= \sum_{w\in W(\t ,\g)}c_w h^{w\la}
$$
for some constants $c_w$. Then we have
$$
\tau_{\la'} (h) \= \sum_{w\in W(\t ,\g)}c_w h^{w\la'}.
$$
To see this, choose a $\la''$ dominating both $\la$ and $\la'$,
then apply the Zuckerman functors $\ph_{\la''}^\la$ and
$\ph_{\la''}^{\la'}$. Proposition 10.44 of \cite{Knapp} gives the
claim. \qed

Write $\pi_\la$ for the discrete series representation with
infinitesimal character $\la$. Let $d(\la):= d_{\pi_\la}$ be the
formal degree then $d(\la)$ is a polynomial in $\la$, more
precisely from \cite{AtSch} we take
$$
d(\la) \= \prod_{\alpha\in\Phi^+(\t ,\g)} \frac{(\alpha ,\la
+\rho)}{(\alpha ,\rho)},
$$
where the ordering $\Phi^+$ is chosen to make $\la$ positive.

Putting things together we see that for $\Re(s)$ small enough the
order of $Z_{P,\tau,\omega}(s)$ at $s=\mu(H_1)$ is
\begin{eqnarray*}
\CO (\mu) &\=& \dim \omega (-1)^{\frac{\dim X}{2}} \chi (X_\Ga)\\
    & {}&\times\sum_{j\in J} \sum_{w\in W(\t ,\g)} \sum_{\pi \in \hat{G}_{ds,j}}
    d(\la_\pi) m_w(\pi_j,\tau ,\mu +w(\la_j-\la_\pi)|_\a),
\end{eqnarray*}
where $\pi_j\in\hat{G}_{ds,j}$ is a fixed element. The function
$\mu \mapsto m_w(\pi_j ,\tau ,\mu)$ takes nonzero values only for
finitely many $\mu$. Since further $\la\mapsto d(\la)$ is a
polynomial it follows that the regularized product
$$
D_{P,\tau ,\omega}(s) := \widehat{\prod_{\mu ,\CO(\mu)\neq 0}}
(s-\mu(H_1))^{\CO(\mu)}
$$
exists. We now have proven the following theorem.

\begin{theorem}
With
$$
\hat{Z}_{P,\tau ,\omega}(s) := Z_{P,\tau,\omega}(s) D_{H,\tau
,\omega}(s)^{-1}
$$
we have
$$
\hat{Z}_{P,\tau ,\omega}(2| \rho_0| -s) \= e^{Q(s)}
\hat{Z}_{P,\tau ,\omega}(s),
$$ where $Q$ is a polynomial.
\qed
\end{theorem}

\begin{proposition} \label{extformel}
Let $\tau$ be a finite dimensional representation of $M$ then the
order of $Z_{P,\tau , \omega}(s)$ at $s=\la(H_1)$ is
$$
(-1)^{\dim(N)} \sum_{\pi \in \hat{G}}N_{\Ga,\omega}(^\theta\pi)
\sum_{q=0}^{\dim(\m\oplus{\n} /\k_M)}(-1)^q \dim(H^q(\m\oplus\n,
K_M ,\pi_K \otimes V_{\breve{\tau}})^\la).
$$
This can also be expressed as
$$
(-1)^{\dim(N)} \sum_{\pi \in \hat{G}}N_{\Ga,\omega}(^\theta\pi)
\sum_{q=0}^{\dim(\m\oplus\n /\k_M)}(-1)^q \dim({\rm Ext}_{(\m
\oplus {\n} ,K_M)}^q (V_{{\tau}} ,V_\pi)^\la).
$$
\end{proposition}

\prf Extend $V_{{\tau}}$ to a $\m \oplus \n$-module by letting
$\n$ act trivially. We then get
$$
H^p(\n,\pi_K) \otimes V_{{\tau}} \cong H^p(\n,\pi_K \otimes
V_{\breve{\tau}}).
$$

The $(\m ,K_M)$-cohomology of the module $H^p(\n,\pi_K \otimes
V_{\tau})$ is the cohomology of the complex $(C^*)$ with
\begin{eqnarray*}
C^q &\=& {\rm Hom}_{K_M}(\wedge^q\p_M ,H^p(\n,\pi_K)\otimes
V_{\breve{\tau}})
\\
        &\=& (\wedge^q\p_M \otimes H^p(\n,\pi_K)\otimes V_{\breve{\tau}})^{K_M},
\end{eqnarray*}
since $\wedge^p\p_M$ is a self-dual $K_M$-module. Therefore we
have an isomorphism of virtual $A$-modules:
$$
\sum_q (-1)^q (H^p(\n,\pi_K)\otimes \wedge^q\p_M \otimes
V_{\breve{\tau}} )^{K_M} \cong \sum_q (-1)^q H^q (\m
,K_M,H^p(\n,\pi_K \otimes V_{\breve{\tau}})).
$$

Now one considers the Hochschild-Serre spectral sequence in the
relative case for the exact sequence of Lie algebras
$$
0 \to \n\to \m \oplus \n\to \m \to 0
$$
and the $(\m\oplus \n,K_M)$-module $\pi \otimes V_{\breve{\tau}}$.
We have
$$
E_2^{p,q} \= H^q(\m ,K_M ,H^p(\n,\pi_K\otimes V_{\breve{\tau}}))
$$
and
$$
E_\infty^{p,q} \= {\rm Gr}^q(H^{p+q}(\m \oplus \n,K_M ,\pi_K
\otimes V_{\breve{\tau}})).
$$
Now the module in question is just
$$
\chi(E_2) \= \sum_{p,q} (-1)^{p+q} E_2^{p,q}.
$$
Since the differentials in the spectral sequence are
$A$-homomorphisms this equals $\chi(E_\infty)$. So we get an
$A$-module isomorphism of virtual $A$-modules
$$
\sum_{p,q} (-1)^{p+q} (H^p(\n,\pi_K)\otimes \wedge^q\p_M \otimes
V_{\breve{\tau}})^{K_M} \cong \sum_j (-1)^j H^j (\m\oplus \n,K_M,
\pi_K \otimes V_{\breve{\tau}}).
$$
The second statement is clear by \cite{BorWall} p.16. \qed

\subsection{The Ruelle zeta function}
The generalized Ruelle zeta function can be described in terms of
the Selberg zeta function as follows.

\begin{theorem}
Let $\Ga$ be  neat and choose a parabolic $P$ of splitrank
one. For $\Re(s)>>0$ define the zeta function
$$
Z_{P,\omega}^R(s) \= \prod_{[\ga]\in {\cal E}_H^p(\Ga)}
\det\left(\begin{array}{c}1-e^{-sl_\ga}\omega(\ga)\end{array}\right)^{\chi_{_1}(X_\ga)},
$$
then $Z_{P,\omega}^R(s)$ extends to a meromorphic function on
$\C$. More precisely, let $\n =\n_\alpha \oplus \n_{2\alpha}$ be
the root space decomposition of $\n$ with respect to the roots of
$(\a ,\g)$ then
$$
Z_{P,\omega}^R(s) = \prod_{q=0}^{\dim \n_\alpha} \prod_{p=0}^{\dim
\n_{2\alpha}}
Z_{P,(\wedge^q\n_\alpha)\otimes(\wedge^p\n_{2\alpha}),\omega}(s+(q+2p)|\alpha|)^{(-1)^{p+q}}.
$$
\end{theorem}

In the case when ${\rm rank}_\R G =1$ this zeta function coincides
with the Ruelle zeta function of the geodesic flow of $X_\Ga$.

\prf For any finite dimensional virtual representation $\tau$ of
$M$ we compute
$$
\log Z_{P,\tau,\omega}(s) = \sum_{[\ga]\in\CE_P^p(\Ga)}
\chi_{_1}(X_\ga)
        \sum_{N\ge 0} \tr(\log(1-e^{-sl_\ga}\ga)|\omega\otimes\tau\otimes S^N(\n))
$$
\begin{eqnarray*}
&=& \sum_{[\ga]\in\CE_P^p(\Ga)} \chi_{_1}(X_\ga)
        \sum_{N\ge 0} \sum_{n\ge 1} \frac 1{n} e^{-sl_\ga n}\tr(\ga |\omega\otimes\tau\otimes S^N(\n))\\
&=& \sum_{[\ga]\in\CE_P(\Ga)} \chi_{_1}(X_\ga)
\frac{e^{-sl_\ga}}{\mu(\ga)}
    \tr\omega(\ga)
    \frac{\tr\tau(b_\ga)}{\det(1-(a_\ga b_\ga)^{-1}|\n)}\\
&=& \sum_{[\ga]\in\CE_P(\Ga)} \chi_{_1}(X_\ga)
\frac{e^{-sl_\ga}}{\mu(\ga)}
    \tr\omega(\ga)
    \frac{\tr\tau(b_\ga)}{\tr((a_\ga b_\ga)^{-1}|\wedge^*\n)}.
\end{eqnarray*}
Since $\n$ is an $M$-module defined over the reals we conclude
that the trace $ \tr((a_\ga m_\ga)^{-1}|\wedge^*\n) $ is a real
number. Therefore it equals its complex conjugate which is
$\tr(a_\ga^{-1} m_\ga|\wedge^*\n)$. Now split into the
contributions from $\n_{\alpha}$ and $\n_{2\alpha}$. The claim now
becomes clear. \qed

\section{The Patterson conjecture}
The content of this section is not strictly needed in the rest of the paper, but it is a result of independent interest.

In this section $G$ continues to be a connected semisimple
Lie group with finite center.
Fix a parabolic $P=MAN$ of splitrank one, ie, $\dim A=1$. Let $\nu \in \a^*$ and let
$({{\sigma}},V_\sigma)$ be a finite dimensional complex
representation of $M$. Consider the principal series
representation $\pi_{\sigma ,\nu}$ on the Hilbert space
$\pi_{\sigma ,\nu}$ of all functions $f$ from $G$ to
$V_\sigma$ such that $f(xman)=a^{-(\nu +\rho)}
\sigma(m)^{-1} f(x)$ and such that the restriction of $f$
to $K$ is an $L^2$-function. Let $\pi_{\sigma,\nu}^\infty$ denote the Fr\'echet space of smooth vectors
and
$\pi_{\sigma,\nu}^{-\infty}$ for its continuous dual. For
$\nu
\in
\a^*$ let
$\bar{\nu}$ denote its complex conjugate with respect to
the real form $\a_0^*$.

We will now formulate the Patterson conjecture \cite{BuOl}.

\begin{theorem} If the order of the Weyl group $W(G,A)$ equals 2, then the cohomology  group $H^p(\Ga ,\pi_{\sigma ,\nu}^{-\infty} \otimes
V_\omega)$ is finite dimensional for every $p\geq 0$ and the vanishing order of the Selberg zeta function
$Z_{P,\tau,\omega}(s+|\rho_0|)$ at $s=\nu(H_1)$ equals
$$
\ord_{s=\nu(H_1)}Z_{P,\tau,\omega}(s+|\rho_0|) \= \chi_{_1}(\Ga ,\pi_{\breve\sigma ,{-\nu}}^{-\infty}
\otimes V_{\omega})
$$
if $\nu \neq 0$ and
$$
\ord_{s=0}Z_{P,\tau,\omega}(s+|\rho_0|)\= \chi_{_1}(\Ga ,\hat{H}_{\breve\sigma,0}^{-\infty}\otimes
V_{\omega}),
$$
where $\hat{H}_{\sigma ,0}^{-\infty}$ is a certain nontrivial
extension of $\pi_{\sigma ,0}^{-\infty}$ with itself.

Further $\chi (\Ga ,\pi_{\sigma ,\nu}^{-\infty}\otimes
V_{\breve{\omega}})$ vanishes.
\end{theorem}

\prf Let $\a_0$ be the real
Lie algebra of $A$ and $\a_0^-$ the negative Weyl chamber with
$\exp(\a_0^-)=A^-$. Let $H_1\in \a_0^-$ be the unique element of
norm $1$. Recall that the vanishing order equals
$$
\sum_{\pi\in\hat{G}} N_{\Ga ,\omega}(\pi) \sum_{q=0}^{\dim
N}\sum_{p=0}^{\dim\p_M} (-1)^{p+q}
\dim(H_q({\n},\pi_K)^{\nu-\rho_0}\otimes \wedge^p\p_M\otimes
V_{\breve{\sigma}})^{K_M}.
$$
Set $\la=\nu-\rho_0$.

\begin{lemma}
If $\la\ne -\rho_0$ then $\a$ acts semisimply on
$H_q({\n},\pi_K)^\la$. The space $H_\rho({\n},\pi_K)^{-\rho_0}$
is annihilated by $(H+\rho_0(H))^2$ for any $H\in\a$.
\end{lemma}

Since $W(G,A)$ is nontrivial, the proof in \cite{BuOl} Prop. 4.1. carries over to our situation.
\qed

To prove the theorem, assume first $\la\ne -\rho_0$. Let
$\pi\in\hat{G}$, then
$$
\sum_{q=0}^{\dim N}\sum_{p=0}^{\dim\p_M} (-1)^{p+q} \dim(
H_q({\n},\pi_K)^\la \otimes \wedge^p\p_M\otimes
V_{\breve{\sigma}} )^{K_M}
$$ $$
=\ \sum_{q=0}^{\dim N} \sum_{p=0}^{\dim\p_M} (-1)^{p+q} \dim(
H^0(\a ,H^p(\m,K_M,H_q({\n},\pi_K) \otimes V_{{\breve{\sigma}}
,-\la}))),
$$
where $V_{{\breve{\sigma}} ,-\la}$ is the representation
space of the representation ${\breve{\sigma}}\otimes
(-\la)$. We want to show that this equals
$$
\sum_{q=0}^\infty \sum_{p=0}^\infty (-1)^{p+q+r}p
\dim H^p(\a\oplus\m,K_M, H_q({\n},\pi_K)\otimes
V_{{\breve{\sigma}} ,-\la}).
$$ 

Since the Hochschild-Serre spectral sequence degenerates for a one dimensional Lie
algebra we get that
$$
\dim H^p(\a\oplus\m,K_M,V)
$$
equals
$$
\dim H^0(\a,H^{p-1}(\m,K_M,V))+\dim H^0(\a,H^{p}(\m,K_M,V)).
$$
This implies
$$
\sum_{p,q\ge 0} (-1)^{p+q+r}p \sum_{b=p-1}^p
\dim H^0(\a,H^b(\m,K_M, H_q({\n},\pi_K)\otimes
V_{{\breve{\sigma}} ,-\la}))
$$ $$
\= \sum_{b,q\ge 0}(-1)^q \sum_{p=b}^{b+1} p\, \dim H^0(\a,H^b(\m,K_M,
H_q({\n},\pi_K)\otimes V_{{\breve{\sigma}} ,-\la})).
$$
It follows that the vanishing order equals
$$
\sum_{q=0}^\infty \sum_{p=0}^\infty (-1)^{p+q+1} p
\dim H^p(\a\oplus\m,K_M, H_q({\n},\pi_K)\otimes
V_{{\breve{\sigma}} ,-\la}).
$$ $$
=\ \sum_{q=0}^\infty \sum_{p=0}^\infty (-1)^{p+q+1} \dim \Ext_{\a\oplus\m,K_M}^p(
H_q({\n},\breve{\pi_K}), V_{{\breve{\sigma}} ,-\la}).
$$

For a $(\g,K)$-module $V$ and a $(\a\oplus\m,M)$-module $U$ we
have \cite{HeSch}:
$$
\Hom_{\g,K}(V,{\rm Ind}_P^G(U)) \ \cong\ \Hom_{\a\oplus\m,K_M}(H_0(\n
,V),U\otimes\C_{\rho_0}),
$$
where $\C_{\rho_0}$ is the one dimensional $A$-module
given by $\rho_0$. 
Thus
$$
\sum_{\pi\in\hat{G}} N_{\Ga,\omega}(\pi) \sum_{p=0}^\infty
\sum_{q=0}^\infty (-1)^{p+q+1} \dim
\Ext_{\a\oplus\m,K_M}^p
(H_q(\n,\breve{\pi_K}),\pi_{{\breve{\sigma}},-\la +\rho_P,K}^\infty)
$$
\begin{eqnarray*}
    &=& \sum_{p=0}^\infty (-1)^{p+1} \dim \Ext_{\g,K}^p (C^\infty(\Ga\bs G,\breve{\omega})_K,
    \pi_{{\breve{\sigma}}, -\la,K}^\infty)
\end{eqnarray*}
Diualizing shows that this equals
$$
\sum_{p\ge0} (-1)^{p+1}\dim \Ext_{\g,K}^p(\pi_{\sigma,\la,K}^\infty,C^\infty(\Ga\bs
G,\omega)_K).
$$
Next by \cite{BorWall}, Chap. I,
$$
\Ext_{\g,K}^p(\pi_{\sigma,\la,K}^\infty,C^\infty(\Ga\bs
G,\omega)_K)\ \cong\ H^p(\g,K,\Hom_\C(\pi_{\sigma,\la,K}^\infty,C^\infty(\Ga\bs G,\omega)_K).
$$
For any two smooth $G$-representations $V,W$ the restriction map gives an isomorphism
$$
\Hom_{ct}(V,W)_K\ \cong\ \Hom_\C(V_K,W_K)_K,
$$
where $\Hom_{ct}$ means continuous homomorphisms. Therefore, using the classical identification of
$(\g,K)$ with differentiable and continuous cohomology as in
\cite{BorWall} we get
\begin{eqnarray*}
\Ext_{\g,K}^q(\pi_{\sigma,\la,K}^\infty,C^\infty(\Ga\bs G,\omega)_K) &\cong &
H^p(\g,K,\Hom_{ct}(\pi_{\sigma,\la}^\infty,C^\infty(\Ga\bs G,\omega))_K)\\
&\cong &
H^p(\g,K,\Hom_{ct}(\pi_{\sigma,\la}^\infty,C^\infty(\Ga\bs G,\omega)))\\
&\cong &
H_d^p(G,\Hom_{ct}(\pi_{\sigma,\la}^\infty,C^\infty(\Ga\bs G,\omega)))\\
&\cong &
H_{ct}^p(G,\Hom_{ct}(\pi_{\sigma,\la}^\infty,C^\infty(\Ga\bs G,\omega)))\\
&\cong &
\Ext_G^p(\pi_{\sigma,\la}^\infty, C^\infty(\Ga\bs G,\omega))\\
&\cong &
\Ext_\Ga^p(\pi_{\sigma,\la}^\infty,\omega)\\
&\cong&
H^p(\Ga,\pi_{\breve\sigma,-\la}^{-\infty}\otimes\omega).
\end{eqnarray*}
This gives the claim for $s\ne 0$. In the case $s=0$ we have to replace $H^0(\a,\cdot)$ by
$H^0(\a^2,\cdot)$, where $\a^2$ means the subalgebra of $U(\a)$ generated by $H^2$, $H\in\a$. Then
the induced representation ${\rm Ind}_P^G(U)$ is replaced by a suitable self-extension.
\qed

\section{Holomorphic torsion}
\subsection{Holomorphic torsion}
Let $D$ be a self-adjoint unbounded operator on a Hilbert space.
Assume that $D$ has eigenvalues $0\le\la_1\le\la_2\le\dots$ and that the zeta function
$$
\zeta_D(s)\=\sum_{j\ge 1: \la_j>0}\la_j^{-s}
$$
converges in some half plane and extends to a meromorphic function on the plane which is regular at $s=0$.
In that case we call $D$ \e{zeta-admissible} and define
$$
\det'(D)\=\exp\(-\zeta_D'(0)\).
$$

Let $E=(E_0 \rightarrow E_1 \rightarrow \dots \rightarrow E_n)$ be
an elliptic complex over a compact smooth manifold $M$. That is,
each $E_j$ is a vector bundle over $M$ and there is differential
operators $d_j:\Ga^\infty(E_j)\to\Ga^\infty(E_{j+1})$ of order one
such that the sequence of principal symbols
$$
\sigma_E(x,\xi)\ :\ 0\to E_{0,x}\to E_{1,x}\to\dots\to E_{n,x}\to
0
$$
is exact whenever $\xi\in T_xM^*$ is nonzero.

Assume each $E_k$ is equipped with a Hermitian metric. Then we can
form the Laplace operators $\Delta_j= d_jd_j^*+d_j^*d_j$ as second
order differential operators. When considered as unbounded
operators on the spaces $L^2(M,E_k)$ these are known to be zeta admissible. Now define the \emph{torsion} of $E$ as $$ \tau_1 (E)
\= \prod_{k=0}^n \det'(\Delta_k)^{k(-1)^{k+1}}.
$$ Note that this definition differs by an exponent 2 from the
original one \cite{RS-RT}.

The corresponding notion in the combinatorial case, i.e. for a
finite CW-complex and the combinatorial Laplacians, was introduced
by Reidemeister in the 1930's who used it to distinguish homotopy
equivalent spaces which are not homeomorphic. The torsion as
defined above, also called {\it analytic torsion}\index{analytic
torsion}, was defined by Ray and Singer in \cite{RS-RT} where they
conjectured that the combinatorial and the analytical torsion
should coincide. This conjecture was later proven independently by
J. Cheeger and W. M\"uller.

For a compact smooth Riemannian manifold $M$ and $E\rightarrow M$
a flat Hermitian vector bundle the complex of $E$-valued forms on
$M$ satisfies the conditions above so that we can define the
torsion $\tau (E)$ via the de Rham complex. Now assume further,
$M$ is K\"{a}hlerian and $E$ holomorphic then we may also consider
the torsion $T(E)$ of the Dolbeault complex $\bar{\partial}\ :\
\Omega^{0,.}(M,E) \rightarrow \Omega^{0,.+1}(M,E)$. This then is
called the {\it holomorphic torsion}\index{holomorphic torsion}.
The holomorphic torsion is of significance in Arakelov theory
\cite{soule} where it is used as normalization factor for families
of Hermitian metrics.

We now define \emph{$L^2$-torsion}. For the following see also
\cite{L}. Let $M$ denote a compact oriented smooth manifold, $\Ga$
its fundamental group and $\tilde{M}$ its universal covering. Let
$E=E_0\rightarrow \dots \rightarrow E_n$ denote an elliptic
complex over $M$ and $\tilde{E}=\tilde{E}_0\rightarrow \dots
\rightarrow \tilde{E}_n$ its pullback to $\tilde{M}$. Assume all
$E_k$ are equipped with Hermitian metrics.

Let $\tilde{\triangle}_p$ and $\triangle_p$ denote the
corresponding Laplacians. The ordinary torsion was defined via the
trace of the complex powers $\triangle_p^s$. The $L^2$-torsion
will instead be defined by considering the complex powers of
$\tilde{\triangle}_p$ and applying a different trace functional.
Write ${\cal F}$ for a fundamental domain of the $\Ga$-action on
$\tilde{M}$ then as a $\Ga$-module we have $$ L^2(\tilde{E}_p)\
\cong\ l^2(\Ga) \otimes L^2(\tilde{E}_p\mid_{\cal F})\ \cong\
l^2(\Ga) \otimes L^2(E_p). $$ The von Neumann algebra $VN(\Ga)$
generated by the right action of $\Ga$ on $l^2(\Ga)$has a
canonical trace making it a type ${\rm II}_1$ von Neumann algebra
if $\Ga$ is infinite \cite{GHJ}. This trace and the canonical
trace on the space $B(L^2(E))$ of bounded linear operators on
$L^2(E)$ define a trace ${\rm tr}_\Ga$ on $VN(\Ga) \otimes
B(L^2(E))$ which makes it a type ${\rm II}_\infty$ von Neumann
algebra. The corresponding dimension function is denoted
$\dim_\Ga$.  Assume for example, a $\Ga$-invariant operator $T$ on
$L^2(E)$ is given as an integral operator with a smooth kernel
$k_T$, then a computation shows
$$
\tr_\Ga(T) \= \int_{\cal F} \tr(k_T(x,x))\ dx.
$$
It follows for the heat operator $e^{-t\tilde{\Delta}_p}$ that $$
\tr_\Ga e^{-t\tilde{\Delta}_p} \= \int_{\cal F}\tr <x\mid
e^{-t\tilde{\Delta}_p}\mid x> dx. $$ From this we read off that
$\tr_\Ga e^{-t\tilde{\Delta}_p}$ satisfies the same small time
asymptotic as $\tr e^{-t\Delta_p}$.

Let $\tilde{\Delta}_p' =
\tilde{\Delta}_p|_{\ker(\tilde{\Delta}_p)^\perp}$. Unfortunately very
little is known about large time asymptotic of
$\tr_\Ga(e^{-t\tilde{\Delta}_p'})$ (see \cite{LL}). Let $$
NS(\Delta_p) \= \sup \{ \alpha \in \R \mid \tr_\Ga e^{-t
\tilde{\Delta}_p'} \= O(t^{-\alpha/2})\ {\rm as}\ t\rightarrow
\infty \} $$ denote the \emph{Novikov-Shubin invariant} of $\Delta_p$
(\cite{GrSh}, \cite{LL}).

Then $NS(\Delta_p)$ is always $\geq 0$; in this section we will \emph{assume} that the Novikov-Shubin invariant of $\Delta_p$ is positive.
This is in general an unproven conjecture. In the cases of our
concern in later sections, however, the operators in question are
homogeneous and it can be proven then that their Novikov-Shubin
invariants are in fact positive. We will consider the integral $$
\zeta_{\Delta_p}^1(s) \= \frac 1{\Ga (s)} \int_0^1 t^{s-1} \tr_\Ga
e^{-t\tilde{\Delta}_p'}\ dt, $$ which converges for $\Re (s) >>0$
and extends to a meromorphic function on the entire plane which is
holomorphic at $s=0$, as is easily shown by using the small time
asymptotic (\cite{BGV},Thm 2.30).

Further the integral
$$
\zeta_{\Delta_p(s)}^2(s) \= \frac 1{\Ga(s)}\int_1^\infty
t^{s-1}\tr_\Ga e^{-t\tilde{\Delta}_p'}\ dt
$$
converges for $\Re (s)<\frac 1{2}NS(\Delta_p)$, so in this region we
define the \emph{$L^2$-zeta function} of $\Delta_p$ as
$$
\zeta_{\Delta_p}^{(2)} (s) \= \zeta_{\Delta_p}^1(s) +
\zeta_{\Delta_p}^2(s).
$$
Assuming the Novikov-Shubin invariant of $\Delta_p$ to be positive
we define the \emph{$L^2$-determinant} of $\Delta_p$ as
$$
{\det}^{(2)}(\Delta_p) \= \exp \left( \left. -\frac{d}{d{s}}
\right|_{s=0} \zeta_{\Delta_p}^{(2)} (s)\right) .
$$

Now let the $L^2$-torsion be defined by $$ T^{(2)}(E) \=
\prod_{p=0}^n {\det}^{(2)}(\Delta_p)^{p(-1)^{p+1}}. $$ Again let $M$
be a K\"ahler manifold and $E\to M$ a flat Hermitian holomorphic
vector bundle then we will write $T^{(2)}_{hol}(E)$ for the
$L^2$-torsion of the Dolbeault complex $\Omega^{0,*}(M,E)$.

\subsection{Computation of Casimir eigenvalues}
In this section we give formulas for Casimir eigenvalues of
irreducible representations of reductive groups which hold under
conditions on the $K$-types.

Let $G$ be a connected semisimple Lie group with finite center.
Fix a maximal compact subgroup $K$ with Cartan involution
$\theta$. Let $P=MAN$ be the Langlands decomposition of a
parabolic subgroup $P$ of $G$. Modulo conjugation we can assume
that $AM$ is stable under $\theta$ and then $K_M=K\cap M$ is a
maximal compact subgroup of $M$. Let $\m =\k_M\oplus\p_M$ be the
corresponding Cartan decomposition of the complexified Lie algebra
of $M$. Let $C_M$ denote the Casimir operator of $M$ induced by
the Killing form on $G$.

\begin{lemma} \label{Caseig1}
Let $(\sigma ,V_\sigma)$ be an irreducible finite dimensional
representation of $M$. Let $(\xi ,V_\xi)$ be an irreducible
unitary representation of $M$ and assume
$$
\sum_{p=0}^{\dim (\p)} (-1)^p \dim (V_\xi \otimes \wedge^p\p_M
\otimes V_{\sigma})^{K_M} \neq 0,
$$
then the Casimir eigenvalues satisfy
 $$
 \xi(C_M) = \sigma(C_M).
 $$
\end{lemma}

\prf Recall that the Killing form of $G$ defines a
$K_M$-isomorphism between $\p_M$ and its dual $\p_M^*$, hence in
the assumption of the lemma we may replace $\p_M$ by $\p_M^*$. Let
$\xi_K$ denote the $(\m ,K_M)$-module of $K_M$-finite vectors in
$V_\xi$ and let $C^q(\xi_K\otimes V_\sigma) =
\Hom_{\k_M}(\wedge^q\p_M ,\xi_K\otimes V_\sigma) =
(\wedge^q\p_M^*\otimes\xi_K\otimes V_\sigma)^{\k_M}$ the standard
complex for the relative Lie algebra cohomology $H^q(\m ,\k_M
,\xi_K\otimes V_\sigma)$. Further
$(\wedge^q\p_M^*\otimes\xi_K\otimes V_\sigma)^{K_M}$ forms the
standard complex for the relative $(\m ,K_M)$-cohomology $H^q(\m
,K_M ,\xi_K\otimes V_\sigma)$. In \cite{BorWall}, p.28 it is shown
that
$$
H^q(\m ,K_M ,\xi_K\otimes V_\sigma)=H^q(\m ,\k_M
,\xi_K\otimes V_\sigma)^{K_M/K_M^0}.
$$
Our assumption
implies $\sum_q(-1)^q \dim H^q(\m ,K_M ,\xi_K\otimes
V_\sigma) \ne 0$, therefore there is a $q$ with $0\ne
H^q(\m ,K_M ,\xi_K\otimes V_\sigma)=H^q(\m ,\k_M
,\xi_K\otimes V_\sigma)^{K_M/K_M^0}$, hence $H^q(\m ,\k_M
,\xi_K\otimes V_\sigma)\ne 0$. Now Proposition 3.1 on page
52 of \cite{BorWall} says that $\pi(C_M)\ne\sigma(C_M)$
implies that $H^q(\m,\k_M,\xi_K\otimes V_\sigma)=0$ for
all $q$. The claim follows. \qed

Let $X=G/K$ be the symmetric space to $G$ and assume that $X$ is
Hermitian\index{Hermitian symmetric space}, i.e., $X$ has a
complex structure which is stable under $G$. Let $\theta$ denote
the Cartan involution fixing $K$ pointwise. Since $X$ is
hermitian it follows that $G$ admits a compact Cartan subgroup
$T\subset K$. We denote the real Lie algebras of $G$, $K$ and $T$
by $\g_0, \k_0$ and $\t_0$ and their complexifications by $\g
,\k$ and $\t$. We will fix a scalar multiple $B$ of the Killing
form. As well, we will
write $B$ for its diagonal, so $B(X) = B(X,X)$. Denote by $\p_0$
the orthocomplement of $\k_0$ in $\g_0$ with respect to $B$; then
via the differential of {\rm exp} the space $\p_0$ is isomorphic
to the real tangent space of $X=G/K$ at the point $eK$. Let $\Phi
(\t ,\g)$ denote the system of roots of $(\t ,\g)$, let $\Phi_c
(\t ,\g) = \Phi (\t ,\k)$ denote the subset of compact roots and
$\Phi_{nc} = \Phi -\Phi_c$ the set of noncompact roots. To any
root $\alpha$ let $\g_\alpha$ denote the corresponding root
space. Fix an ordering $\Phi^+$ on $ \Phi = \Phi (\t ,\g)$ and
let $\p_\pm = \bigoplus_{\alpha \in \Phi_{nc}^+}\g_{\pm \alpha}$.
Then the complexification $\p$ of $\p_0$ splits as $\p = \p_+
\oplus \p_-$ and the ordering can be chosen such that this
decomposition corresponds via {\rm exp} to the decomposition of
the complexified tangent space of $X$ into holomorphic and
antiholomorphic part. By Lemma 2.2.3 of \cite{lefschetz} we can,
replacing $G$ by a double cover if necessary, assume that the
adjoint homomorphism $K\to SO(\p)$ factors over ${\rm Spin}(\p)$.

\begin{lemma} \label{Caseig2}
Let $(\tau ,V_\tau)$ denote an irreducible representation of $K$.
Assume $X$ Hermitian and let $(\pi ,W_\pi)$ be an irreducible
unitary representation of $G$ and assume that
$$
\sum_{p=0}^{\dim \p_-} (-1)^p \dim(W_\pi \otimes
\wedge^p\p_-\otimes V_{{\tau}})^K \neq 0
$$
then we have
$$
\pi(C) \= {\tau}\otimes\epsilon(C_K) -B(\rho)+B(\rho_K),
$$
where $\epsilon$ is the one dimensional representation of $K$
satisfying $\epsilon\otimes\epsilon\cong\wedge^{top}\p_+$.
\end{lemma}

\prf This follows from Lemma 2.1.1 of \cite{lefschetz}. \qed

\subsection{The local trace of the heat kernel}
On a Hermitian globally symmetric space the heat operator is given
by convolution with a function on the group of isometries. In this
section we determine the trace of this function on any irreducible
unitary representation.

Let $H$ denote a $\theta$-stable Cartan subgroup of $G$ so $H=A B$
where $A$ is the connected split component and $B$ compact. The
dimension of $A$ is called the split rank of H. Let $\a$ denote
the complex Lie algebra of $A$. Then $\a$ is an abelian subspace
of $\p=\p_+ \oplus \p_-$. Let $X\mapsto X^c$ denote the complex
conjugation on $\g$ according to the real form $\g_0$. The next
lemma shows that $\a$ lies skew to the decomposition $\p =\p_+
\oplus \p_-$.

\begin{lemma} Let $Pr_\pm$ denote the projections from $\p$ to $\p_\pm$.
Then we have $\dim Pr_+(\a)$ = {\rm dim} $Pr_-(\a)$ = {\rm dim}
$\a$, or, what amounts to the same: $\a \cap \p_\pm = 0$.
\end{lemma}

Proof: $\a$ consists of semisimple elements whereas $\p_\pm$
consists of nilpotent elements only.
 \qed

Now let $\a'$ denote the orthocomplement of $\a$ in $Pr_+(\a)
\oplus Pr_-(\a)$. For later use we write ${\cal V} = \a \oplus \a'
= {\cal V}_+ \oplus {\cal V}_-$, where ${\cal V}_\pm = {\cal
V}\cap \p_\pm = Pr_\pm (\a)$.

Let $(\tau,V_\tau)$ be a finite dimensional representation of the
compact group $K$ and let $E_\tau = G\times_K V_\tau = (G\times
V)/K$ be the corresponding $G$-homogeneous vector bundle over $X$.
The sections of the bundle $E_\tau$ are as a $G$-module given by
the space of $K$-invariants
 $$
 \Ga^\infty(E_\tau) \= \left( C^\infty(G)\otimes V_\tau\right)^K,
 $$
where $K$ acts on $C^\infty(G)\otimes V_\tau$ by $k.(f(x)\otimes
v) = f(xk^{-1})\otimes \tau(k)v$. and $G$ acts on $\left(
C^\infty(G)\otimes V_\tau\right)^K$ by left translations on the
first factor.

Consider the convolution algebra $C_c^\infty(G)$ of compactly
supported smooth functions on $G$. Let $K\times K$ act on it by
right and left translations and on $\End_\C(V_\tau)$ by
$(k_1,k_2)T=\tau(k_1)T\tau(k_2^{-1})$. Then the space of
invariants
 $$
 \left(C_c^\infty(G)\otimes \End_\C(V_\tau)\right)^{K\times K}
 $$
is seen to be an algebra again and it acts on
$\Ga^\infty(E_\tau)=\left( C^\infty(G)\otimes V_\tau\right)^K$ by
convolution on the first factor and in the obvious way on the
second. This describes all $G$-invariant smoothing operators on
$E_\tau$ which extend to all sections. $G$-invariant smoothing
operators which not extend to all sections will be given by
Schwartz kernels in
 $$
 \left(C^\infty(G)\otimes \End_\C(V_\tau)\right)^{K\times K}.
 $$

Let $n=2m$ denote the real dimension of $X$ and for $0\leq p,q
\leq m$ let $\Omega^{p,q}(X)$ denote the space of smooth
$(p,q)$-forms on $X$. The above calculus holds for the space of
sections $\Omega^{p,q}(X)$. D.Barbasch and H. Moscovici have shown
in \cite{BM} that the heat operator $e^{-t\triangle_{p,q}}$ has a
smooth kernel $h_t^{p,q}$ of rapid decay in
 $$
 ( C^\infty (G)\otimes {\rm End}
 (\wedge^p\p_+ \otimes \wedge^q\p_-))^{K\times K}.
 $$
Now fix p and set for $t>0$ $$ f_t^p \= \sum_{q=0}^m q(-1)^{q+1}\
{\rm tr}\ h_t^{p,q}, $$ where tr means the trace in ${\rm End}
(\wedge^p\p_+ \otimes \wedge^q\p_-)$.

We want to compute the trace of $f_t^p$ on the principal series
representations. To this end let $H=A B$ as above and let $P$
denote a parabolic subgroup of $G$ with Langlands decomposition
$P=MA N$. Let $(\xi ,W_\xi)$ denote an irreducible unitary
representation of $M$, $e^\nu$ a quasicharacter of $A$ and set
$\pi_{\xi ,\nu} = {\rm Ind}_P^G (\xi \otimes e^{\nu + \rho_{_P}}\otimes
1)$, where $\rho_{_P}$ is the half of the sum of the $P$-positive
roots.

Let $C$ denote the Casimir operator of $G$ attached to the form
$B$.

\begin{proposition} \label{trace1}
The trace of $f_t^p$ under $\pi_{\xi ,\nu}$ vanishes if $\dim \a
> 1 $. If ${\rm  dim\ } \a =1$ it equals
$$
e^{t\pi_{\xi ,\nu}(C)} \sum_{q=0}^{\dim(\p_-)-1} (-1)^q\
\dim\left( W_\xi \otimes \wedge^p\p_+ \otimes \wedge^q(\a^\perp
\cap \p_-)\right)^{K\cap M},
$$
where $\a^\perp$ is the orthocomplement of $\a$ in $\p$.
\end{proposition}

Proof: As before let $\a'$ be the span of all $X-X^c$, where
$X+X^c\in \a,\ X\in \p_+$ then ${\cal V} = \a \oplus \a' = \CV_+
\oplus \CV_-$ where $\CV_\pm = \CV \cap \p_\pm$. The group
$K_M=K\cap M$ acts trivially on $\a$ so for $x\in K_M$ we have
$X+X^c = {\rm Ad}(x)(X+X^c) = {\rm Ad}(x)X + {\rm Ad}(x)X^c$.
Since $K$ respects the decomposition $\p=\p_-\oplus \p_+$, we
conclude that $K_M$ acts trivially on $\CV$, hence on $\CV_-$.
Let $r=\dim\a=\dim\CV_-$. As a $K_M$-module we have $$
\begin{array}{cll}
\wedge^p\p_-    &       \=       &       \displaystyle
                                        \sum_{a+b=q} \wedge^b \CV_- \otimes
                                        \wedge^a \CV_-^\perp\\
                &       \=       &       \displaystyle
                                        \sum_{a+b=q} \binom{r}{b}
                                        \wedge^a \CV_-^\perp,
\end{array}
$$
where $\CV_-^\perp = \a^\perp \cap \p_-$. By definition we get
$$
\begin{array}{cll}
\tr\pi_{\xi ,\nu}(f_t^p)    & \= &   \displaystyle
                                \tr \pi_{\xi ,\nu}\left(\sum_{q=0}^m q(-1)^{q+1}
                                        h_t^{p,q}\right)\\
                & \= &   \displaystyle
                        e^{t\pi_{\xi ,\nu}(C)} \sum_{q=0}^m q(-1)^{q+1}
                                \dim(V_{\pi_{\xi ,\nu}} \otimes
                                \wedge^p\p_+ \otimes \wedge^q\p_-)^K
\end{array} $$
By Frobenius reciprocity this equals
\begin{eqnarray*}
                & {} &   \displaystyle
                        e^{t\pi_{\xi ,\nu}(C)} \sum_{q=0}^m q(-1)^{q+1}
                                \dim(W_{\pi_{\xi}} \otimes
                                \wedge^p\p_+ \otimes \wedge^q\p_-)^{K\cap M}\\
                & \= &   \displaystyle
                        e^{t\pi_{\xi ,\nu}(C)} \sum_{q=0}^m \sum_{a=0}^q q(-1)^{q+1}
                                \binom{r}{q-a} \dim(W_{\pi_{\xi}} \otimes
                                \wedge^p\p_+ \otimes
                                \wedge^a\CV_-^\perp)^{K\cap M}\\
                & \= &   \displaystyle
                        e^{t\pi_{\xi ,\nu}(C)} \sum_{a=0}^m \sum_{q=a}^m q(-1)^{q+1}
                                \binom{r}{q-a} \dim(W_{\pi_{\xi}} \otimes
                                \wedge^p\p_+ \otimes
                                \wedge^a\CV_-^\perp)^{K\cap M}
\end{eqnarray*}
By taking into account $a\leq m-r$ we get
 $$
 \sum_{q=a}^m q(-1)^q \binom{r}{q-a} \=
 \left\{ \begin{array}{cl}(-1)^{a+1}&{\rm if}\ r=1,\\ 0 & {\rm if}\ r>1,\end{array}\right.
 $$
 and
the claim follows. \qed

We will now combine this result with the computation of Casimir
eigenvalues. This requires an analysis whether the homogeneous
space $M/K_M$ embedded into $X=G/K$ inherits the complex structure
or not.

Fix a $\theta$-stable Cartan subgroup $H=A B$ with $\dim(A) =1$
and a parabolic $P=MA N$. Fix a system of positive roots
$\Phi^+=\Phi^+(\g,\h)$ in $\Phi(\g,\h)$ such that for
$\alpha\in\Phi^+$ and $\alpha$ nonimaginary it follows
$\alpha^c\in\Phi^+$. Further assume that $\Phi^+$ is compatible
with the choice of $P$, i.e., for any $\alpha\in\Phi^+$ the
restriction $\alpha|_{\a}$ is either zero or positive. Let $\rho$
denote the half sum of positive roots. For $\xi \in \hat{M}$ let
$\la_\xi \in \b^*$ denote the infinitesimal character of $\xi$.
Recall that we have
$$ \pi_{\xi ,\nu}(C) \= B(\nu)+B(\la_\xi)-B(\rho). $$

\begin{lemma}
There exists a unique real root $\alpha_r\in\Phi^+$.
\end{lemma}

\prf Recall that a root $\alpha\in\Phi(\g,\h)$ is real if and only
if it annihilates $\b=Lie_\C(B)$. Hence the real roots are
elements of $\a^*$ which is one dimensional, so, if there were
two positive real roots, one would be positive multiple of the
other which is absurd.

To the existence. We have that $\dim X=\dim\a +\dim\n$ is even,
hence $\dim\n$ is odd. Now $\n=\oplus_\alpha\g_\alpha$, where the
sum runs over all $\alpha\in\Phi^+$ with $\alpha|_\a\ne 0$. The
complex conjugation permutes the $\g_\alpha$ and for any nonreal
$\alpha$ we have $(\g_\alpha)^c\ne\g_\alpha$, hence the nonreal
roots pair up, thus $\n$ can only be odd dimensional if there is a
real root.\qed

Let $c=c(H)$ denote the number of positive restricted roots in
$\Phi(\g,\a)$. Let $\Phi^+(\g,\a)$ denote the subset of positive
restricted roots.

\begin{lemma}
There are three possibilities:
\begin{itemize}
\item
$c=1$ and $\Phi^+(\a,\g)=\{ \alpha_r\}$,
\item
$c=2$ and $\Phi^+(\a,\g)=\{ \frac 1{2}\alpha_r, \alpha_r\}$ and
\item
$c=3$ and $\Phi^+(\a,\g)=\{ \frac 1{2}\alpha_r, \alpha_r,
\frac{3}{2}\alpha_r\}$.
\end{itemize}
\end{lemma}

\prf For any root $\alpha\in\Phi(\h ,\g)$ we have that $2B(\alpha
,\alpha_r)/B(\alpha_r)$ can only take the values $0,1,2,3$, so
the only possible restricted roots in $\Phi^+(\a,\g)$ are
$\frac 1{2}\alpha_r, \alpha_r, \frac{3}{2}\alpha_r$.

If $\frac{3}{2}\alpha_r\in \Phi^+(\a,\g)$ then there is a root
$\beta$ in $\Phi(\h,\g)$ with $\beta |_\a
=\frac{3}{2}\alpha_r|_\a$. Then $B(\alpha_r ,\beta)>0$ hence
$\eta =\beta-\alpha_r$ is a root. Since $\eta |_\a
=\frac 1{2}\alpha_r$ we get that then $\frac 1{2}\alpha_r\in\Phi(\a
,\g)$. From this the claim follows.
\qed

Recall that a {\it central isogeny}\index{central isogeny} $\ph
:L_1\to L_2$ of Lie groups is a surjective homomorphism with
finite kernel which lies in the center of $L_1$.

\begin{lemma}
There is a central isogeny $M_1\times M_2\to M$ such that the
inverse image of $K_M$ is of the form $K_{M_1}\times K_{M_2}$
where $K_{M_j}$ is a maximal compact subgroup of $M_j$ for
$j=1,2$ and such that $\p_M = \p_{M_1}\oplus\p_{M_2}$ as
$K_{M_1}\times K_{M_2} $-module and
\begin{itemize}
\item
the map $X\mapsto [X,X_{\alpha_r}]$ induces a $K_M$-isomorphism
$$
\p_{M_1}\ \cong\ [\p_M ,\g_{\alpha_r}],
$$
\item
with $\p_{M_2 ,\pm}:=\p_{M_2}\cap\p_\pm$ we have
$$
\p_{M_2} =\p_{M_2,+}\oplus\p_{M_2,-}.
$$
\end{itemize}
The latter point implies that the symmetric space $M_2/K_{M_2}$,
naturally embedded into $G/K$, inherits the complex structure.
This in particular implies that $M_2$ is orientation preserving.
Further we have $\p_{M_1}\cap\p_+ =\p_{M_1}\cap\p_- = 0$, which
in turn implies that the symmetric space $M_1/K_{M_1}$ does not
inherit the complex structure of $G/K$.
\end{lemma}

\prf At first we reduce the proof to the case that the center $Z$
of $G$ is trivial. So assume the proposition proved for $G/Z$ then
the covering $M_1\times M_2\to M$ is gotten by pullback from that
of $M/Z$. Thus we may assume that $G$ has trivial center.

Let $H$ be a generator of $\a_0$. Write $H=Y+Y^c$ for some $Y\in
\p_+$. According to the root space decomposition $\g = \a \oplus
\k_M \oplus \p_M \oplus \n \oplus \theta(\n)$ we write
$Y=Y_a+Y_k+Y_p+Y_n+Y_{\theta(n)}$. Because of $\theta(Y)=-Y$ it
follows $Y_k=0$ and $Y_{\theta(n)}=-\theta(Y_n)$. For arbitrary
$k\in K_M$ we have $\Ad(k)Y=Y$ since $\Ad(k)H=H$ and the
projection $Pr_+$ is $K_M$-equivariant. Since the root space
decomposition is stable under $K_M$ it follows $\Ad(K)Y_*=Y_*$ for
$*=a,p,n$. The group $M$ has a compact Cartan, hence $K_M$ cannot
act trivially on any nontrivial element of $\p_M$, so $Y_p=0$ and
$Y_n\in \g_{\al_r}$. 
Since $Y=Y_a+Y_n-\theta(Y_n)$ and $Y\notin\a$
it follows that $Y_n\ne 0$, so $Y_n$ generates $\g_{\alpha_r}$ and
so $K_M$ acts trivially on $\g_{\al_r}$.

Let $\p_{M_2}\subset\p_M$ by definition be the kernel of the map
$X\mapsto [X,Y_n]$. Let $\p_{M_1}$ be its orthocomplement in
$\p_M$. The group $K_M$ stabilizes $Y_n$, so it follows that $K_M$
leaves $\p_{M_2}$ stable hence the orthogonal decomposition $\p_M
=\p_{M_1}\oplus\p_{M_2}$ is $K_M$-stable. Thus the symmetric space
$M/K_M$ decomposes into a product accordingly and so does the
image of $M$ in the group of isometries of $M/K_M$. It follows
that the Lie algebra $\m$ of $M$ splits as a direct sum of ideals
$\m=\m_0\oplus\m_1\oplus\m_2$, where $\m_0$ is the Lie algebra of
the kernel $M_0$ of the map $M\to{\rm
Iso}(M/K_M)=\overline{M_1}\times\overline{M_2}$. Let $L=AM$ be the
Levi component and let $L_c$ denote the centralizer of $\a$ in the
group ${\rm Int}(\g)$. Then $L_c$ is a connected complex group
and, since $G$ has trivial center, $L$ injects into $L_c$. The Lie
algebra $\l$ of $L$ resp. $L_c$ decomposes as
$l=\a\oplus\m_0\oplus\m_1\oplus\m_2$ and so there is an isogeny
$A_c\times M_{0,c}\times M_{1,c}\times M_{2,c}\to L_c$. The group
$A$ injects into $A_c$ and we have $A_c=A\times U$ where $U$ is a
compact one dimensional torus. Let $\hat{M_1}\subset M$ be the
preimage of $\overline{M_1}$ and define $\hat{M_2}$ analogously.
Then $\hat{M_1}$ and $\hat{M_2}$ are subgroups of $L_c$ and we may
define $M_1$ to be the preimage of $\hat{M_1}$ in $U\times
M_{0,c}\times M_{1,c}$ and $M_2$ to be the preimage of $\hat{M_2}$
in $M_{2,c}$. Then the map $M_1\times M_2\to M$ is a central
isogeny.

By definition we already have $\p_{M_1}\cong [\p_M
,\g_{\alpha_r}]$. Now let $X\in\p_{M_2}$, then $[X,Y]
=[X,Y_n]+\theta([X,Y_n]) =0$ since $Y_n\in\g_{\alpha_r}$. By
$[X,H]=0$ this implies $[X,Y^c]=0$ and so
$[Pr_\pm(X),H]=0=[Pr_\pm(X),Y]=[Pr_\pm(X),Y^c]$ for $Pr_\pm$
denoting the projection $\p\to\p_\pm$. Hence
$Pr_\pm(X)\in\p_{M_2}$ which implies
$$
\p_{M_2}\= \p_{M_2 ,+}\oplus\p_{M_2 ,-}.
$$
The lemma is proven.
\qed

Let
$$
\n_{\al_r} := \bigoplus_{{\alpha\in\Phi(\h ,\g)}\ {\alpha |_\a =\al_r
|_\a}}\g_{\al_r}
$$
and define $\n_{\frac 1{2}\al_r}$ and $\n_{\frac{3}{2}\al_r}$
analogously. Let $\n^\al_r := \n_{\frac 1{2}\al_r}\oplus
\n_{\frac{3}{2}\al_r}$, then $\n=\n_{\al_r}\oplus\n^\al_r$.

\begin{proposition}\label{zerlp}
There is a $K_M$-stable subspace $\n_-$ of $\n^{\al_r}$ such that
as $K_M$-module $\n^\al_r \cong \n_-\oplus\tilde{\n_-}$, where
$\tilde{.}$ denotes the contragredient, and as $K_M$-module:
$$
\p_-\ \cong\ \C \oplus \p_{M_1}\oplus\p_{M_2 ,-}\oplus\n_- .
$$
\end{proposition}

\prf The map $a+p+n\mapsto a+p+n-\theta(n)$ induces a
$K_M$-isomorphism $$ \a\oplus\p_M\oplus\n \ \cong\ \p . $$ We now
prove that there is a $K_M$-isomorphism $\n_{\al_r} \cong \g_{\al_r}\oplus
[\p_M ,\g_{\al_r}]$. For this let $\alpha \ne\al_r$ be any root in
$\Phi(\h,\g)$ such that $\alpha |_\a =\al_r |_\a$. Then $B(\alpha
,\al_r)>0$ and hence $\beta :=\alpha -\al_r$ is a root. The root
$\beta$ is imaginary. Assume $\beta$ is compact, then
$\g_\beta\subset\k_M$ and we have $[\g_{\al_r} ,\g_\beta ]=\g_{\al_r}$
which contradicts $[\k_M ,\g_{\al_r} ]=0$. It follows that $\beta$ is
noncompact and so $\g_\beta\subset\p_M$, which implies $\n_{\al_r}
=\g_{\al_r}\oplus [\p_M,\g_{\al_r} ]$.

Now the last lemma implies the assertion.
\qed

We have $M\leftarrow M_1\times M_2$ and so any irreducible
representation $\xi$ of $M$ can be pulled back to $M_1\times M_2$
and be written as a tensor product $\xi =\xi_1\otimes\xi_2$. Let
$\xi$ be irreducible admissible then Proposition \ref{trace1} and
Proposition \ref{zerlp} imply
$$
\tr\pi_{\xi ,\nu}(f_t^0) \= e^{t\pi_{\xi ,\nu}(C)}\sum_{q\ge 0}
(-1)^q \dim\left(
W_\xi\otimes\wedge^q(\p_{M_1}\oplus\p_{M_2,-}\oplus\n_-)\right)^{K_M}.
$$
For $c\ge 0$ let
$$
\wedge^c\n_- \= \bigoplus_{i\in I_c}\sigma_{i}^c\otimes\tau_{i}^c
$$\index{$\sigma_i^c$}\index{$\tau_i^c$}
be a decomposition as $K_M\leftarrow K_{M_1}\times
K_{M_2}$-module of $\n_-$ where $\tau_i^c$ is an irreducible
$K_{M_2}$-module and $\sigma_i^c$ is the image of the projection
$\n^\al_r =\n_-\oplus\tilde{\n_-}\to\n_-$ of an irreducible
$M_1$-submodule of $\n^\al_r$. Then we conclude $$ \tr\pi_{\xi
,\nu}(f_t^0) \= e^{t\pi_{\xi ,\nu}(C)}\sum_{a,b,c\ge 0}
\sum_{i\in I_c} (-1)^{a+b+c} \dim\left( W_{\xi_1}\otimes
\wedge^a\p_{M_1}\otimes\sigma_{i}^c\right)^{K_{M_1}} $$ $$
\hspace{190pt}\times\dim\left( W_{\xi_2}\otimes
\wedge^b\p_{M_2,-}\otimes\tau_{i}^c\right)^{K_{M_2}}. $$

Fix an irreducible representation $\xi =\xi_1\otimes\xi_2$ of $M$.
Let $\xi':=\tilde{\xi_1}\otimes\xi_2$, where $\tilde{\xi_1}$ is
the contragredient representation to $\xi_1$.

We will use the following notation: For an irreducible
representation $\pi$ we denote its infinitesimal character by
$\wedge_\pi$. We will identify $\wedge_\pi$ to a corresponding
element in the dual of a Cartan subalgebra (modulo the Weyl
group), so that it makes sense to write an expression like
$B(\wedge_\pi)$. In the case of $\sigma_i^c$ we have the
following situation: either $\sigma_i^c$ already extends to a
representation of $M_1$ or $\sigma_i^c\oplus\tilde{\sigma_i^c}$
does. In either case we write $\wedge_{\sigma_i^c}$ for the
corresponding infinitesimal character with respect to the group
$M_1$.

\begin{lemma}
If $\xi\cong\xi'$ and $\tr\pi_{\xi ,\nu}(f_t^0)\ne 0$ for some
$t>0$ then
$$
\pi_{\xi ,\nu}(C) \=
B(\wedge_{\sigma_i^c})+B(\wedge_{\tau_{i}^{c}\otimes\epsilon})+B(\nu)-B(\rho)
$$
for some $c\ge 0$ and some $i\in I_c$. Here $\epsilon$ is the one
dimensional representation of $K_{M_2}$ such that
$\epsilon\otimes\epsilon\cong \wedge^{top}\p_{M_2}$.
\end{lemma}

\prf Generally we have
\begin{eqnarray*}
\pi_{\xi ,\nu}(C) &=& B(\wedge_{\pi_{\xi ,\nu}})-B(\rho)\\
&=& B(\wedge_\xi)+B(\nu) -B(\rho)\\
&=& B(\wedge_{\xi_1})+B(\wedge_{\xi_2}) +B(\nu) -B(\rho).
\end{eqnarray*}

Now $\tr\pi_{\xi ,\nu}(f_t^0)\ne 0$ implies firstly
$$
\sum_{a\ge 0} (-1)^a \dim\left( W_{\xi_1}\otimes
\wedge^a\p_{M_1}\otimes\sigma_i^c\right)^{K_{M_1}}\ne 0
$$
for some $c,i$. Now either $\sigma_i^c$ already extends to an
irreducible representation of $M_1$ or
$\sigma_i^c\oplus\tilde{\sigma_i^c}$ does. In the second case we
get
$$
\sum_{a\ge 0} (-1)^a \dim\left( W_{\xi_1}\otimes
\wedge^a\p_{M_1}\otimes(\sigma_i^c\oplus\tilde{\sigma_i^c})\right)^{K_{M_1}}
$$ $$
= 2\sum_{a\ge 0} (-1)^a \dim\left( W_{\xi_1}\otimes
\wedge^a\p_{M_1}\otimes\sigma_i^c\right)^{K_{M_1}}\ne 0.
$$
So in either case Lemma \ref{Caseig1} implies
$$
B(\wedge_{\xi_1}) -B(\rho_{M_1}) \= \xi_1(C_M) \= \sigma_i^c(C_M)
\= B(\wedge_{\sigma_i^c})-B(\rho_M),
$$
hence $B(\wedge_{\xi_1})=B(\wedge_{\sigma_1^c})$.

Next $\tr\pi_{\xi ,\nu}(f_t^0)\ne 0$ implies
$$
\sum_{b\ge 0}(-1)^b\dim\left(
W_{\xi_2}\otimes\wedge^b\p_{M_2,-}\otimes\tau_{i}^{c}\right)^{K_{M_2}}\ne
0
$$
for the same $c,i$. In that case with $\tau =\tau_{i}^{c}$ Lemma
\ref{Caseig2} implies
$$
\xi_2(C_{M_2}) = \tau\otimes\epsilon(C_{K_{M_2}})-B(\rho_{M_2})
+B(\rho_{K_{M_2}}).
$$
Since $\xi_2(C_{M_2})=B(\wedge_{\xi_2})-B(\rho_{M_2})$ we conclude
$$
B(\wedge_{\xi_2}) \=
\tau\otimes\epsilon(C_{K_{M_2}})+B(\rho_{K_{M_2}}) \=
B(\wedge_{\tau\otimes\epsilon}).
$$
\qed

Along the same lines we get

\begin{lemma}
If $\xi$ is not isomorphic to $\xi'$ and $\tr\pi_{\xi ,\nu}(f_t^0)
+\tr\pi_{\xi' ,\nu}(f_t^0)\ne 0$ for some $t>0$ then
$$
\pi_{\xi ,\nu}(C) \=
B(\wedge_{\sigma_i^c})+B(\wedge_{\tau_{i}^{c}\otimes\epsilon})+B(\nu)-B(\rho)
$$
for some $c\ge 0$ and some $i\in I_c$.
\end{lemma}

We abbreviate $s_i^c := B(\wedge_{\sigma_i^c})+
B(\wedge_{\tau_i^c\otimes\epsilon})$.\index{$s_i^c$} We have shown
that if $\xi\cong\xi'$ then
 $$
 \tr\pi_{\xi ,\nu}(f_t^0)
 $$
equals
 $$ e^{t(B(\nu)-B(\rho))}\sum_{a,b,c\ge 0} \sum_{i\in I_c}
e^{t(s_i^c)}(-1)^{a+b+c} \dim\left( W_{\xi_1}\otimes
\wedge^a\p_{M_1}\otimes\sigma_{i}^c\right)^{K_{M_1}}
 $$ $$
\hspace{190pt}\times\dim\left( W_{\xi_2}\otimes
\wedge^b\p_{M_2,-}\otimes\tau_{i}^c\right)^{K_{M_2}}.
 $$
Define the Fourier transform of $f_t^0$ by
 $$
 \hat{f_t^0}_H(\nu,b^*)\= \tr \pi_{\xi_{b^*},\nu}(f_t^0).
 $$
Let $B^*$ be the character group of $B$. According to $M=M_1\times
M_2$ we can write $B=B_1\times B_2$ and so we see that any
character $b^*$ of $B$ decomposes as $b_1^*\times b_2^*$. Let
${b^*}':=\overline{b^*_1}\times b_2^*$. We get

\begin{lemma}
The sum of Fourier transforms
 $$
\hat{f_t^0}_H(\nu ,b^*)+\hat{f_t^0}_H(\nu ,{b^*}')
 $$
equals
 $$
 e^{t(B(\nu)-B(\rho))}\sum_{a,b,c\ge 0} \sum_{i\in I_c}
e^{ts_i^c}(-1)^{a+b+c} \dim\left( (V_{b_1^*}\oplus
V_{\overline{b_1^*}})\otimes
\wedge^a\p_{M_1}\otimes\sigma_{i}^c\right)^{K_{M_1}}
 $$ $$
\hspace{190pt}\times\dim\left( V_{b_2^*}\otimes
\wedge^b\p_{M_2,-}\otimes\tau_{i}^c\right)^{K_{M_2}}.
 $$
\end{lemma}

\subsection{The global trace of the heat kernel}

Fix a discrete cocompact torsion-free subgroup $\Ga$ of $G$ then
the quotient $X_\Ga :=\Ga\bs X = \Ga\bs G/K$ is a non-euclidean
Hermitian locally symmetric space.

Since $X$ is contractible and $\Ga$ acts freely on $X$ the group
$\Ga$ equals the fundamental group of $X_\Ga$. Let $(\omega
,V_\omega)$ denote a finite dimensional unitary representation of
$\Ga$.

D. Barbasch and H. Moscovici have shown in \cite{BM} that the
function $f_t^0$ satisfies the conditions to be plugged into the
trace formula. So in order to compute $\tr R_{\Ga ,\omega}(f_t^0)$
we have to compute the orbital integrals of $f_t^0$. At first let
$h\in G$ be a nonelliptic semisimple element. Since the trace of
$f_t^0$ vanishes on principal series representations which do not
come from splitrank one Cartan subgroups, we see that
$\CO_h(f_t^0)=0$ unless $h\in H$, where $H$ is a splitrank one
Cartan. Write $H=AB$ as before. We have $h=a_hb_h$ and since $h$
is nonelliptic it follows that $a_h$ is regular in $A$. We say
that $h$ is \emph{split regular} in $H$. We can choose a parabolic
$P=MAN$ such that $a_h$ lies in the negative Weyl chamber $A^-$.

Let $V$ denote a finite dimensional complex vector space and let
$A$ be an endomorphism of $V$. Let $\det(A)$ denote the
determinant of $A$, which is the product over all eigenvalues of
$A$ with algebraic multiplicities. Let $\det'(A)$ be the product
of all nonzero eigenvalues with algebraic multiplicities.

In \cite{HC-HA1}, sec 17 Harish-Chandra has shown that for $h_0\in
H$
 $$
\CO_{h_0}(f_t^0) \= \frac{\varpi_{h_0} (\
'F_{f_t^0}^H(h))|_{h=h_0}}
                {c_{h_0} h_0^{\rho_P} \det'(1-h_0^{-1}|(\g /\h)^+)},
 $$
where $(\g /\h)^+$ is the positive part of the root space
decomposition of a compatible ordering. Further $\varpi_h$ is the
differential operator attached to $h$ as follows. Let $\g_h$
denote the centralizer of $h$ in $\g$ and let $\Phi^+(\g_h,\h)$
the positive roots then
 $$
 \varpi_h\= \prod_{\alpha \in\Phi^+(\g_h,\h)} H_\alpha,
 $$
where $H_\alpha$ is the element in $\h$ dual to $\alpha$ via the
bilinear form $B$. In comparison to other sources the formula
above for the orbital integral lacks a factor $[G_h:G_h^0]$ which
doesn't occur because of the choices of Haar measures made.  We
assume the ordering to come from an ordering of $\Phi(\b,\m)$
which is such that for the root space decomposition $\p_M
=\p_M^+\oplus \p_M^-$ it holds $\p_{M_2}\cap\p_M^+ =\p_{M_2,-}$.
For short we will henceforth write $\varpi_{h_0}F(h_0)$ instead of
$\varpi_{h_0}F(h)|_{h=h_0}$ for any function $F$.

Our results on the Fourier transform of $f_t^0$ together with a
computation as in the proof of Lemma 2.3.6 of \cite{lefschetz} imply that
$'F_{f_t^0}^H(h)$ equals
$$
\frac{e^{-l_h^2/4t}}{\sqrt{4\pi t}}
 e^{-tB(\rho)}
\det(1-h^{-1} | \k_M^+\oplus \p_{M_1}^+) \sum_{c\ge 0,\ i\in I_c}
(-1)^c e^{ts_i^c} \tr(b_h^{-1} | \sigma_i^c\oplus \tau_i^c).
$$
Note that this computation involves a summation over $B^*$ and
thus we may replace $\sigma_i^c$ by its dual $\breve{\sigma_i^c}$
without changing the result.

Now let $(\tau ,V_\tau)$ be a finite dimensional unitary
representation of $K_M$ and define for $b\in B$ the monodromy
factor: $$ L^M(b,\tau) \= \frac{\varpi_b( \det(1-b|(\k_M/|\b)^+)
\tr(\tau(b)))}
                {\varpi_b( \det(1-b|(\m/\b)^+))}.
$$

Note that the expression $\varpi_b( \det (1-b|(\m /\b)^+)) $
equals
$$
|W(\m_b ,\b)|\prod_{\alpha \in \Phi_b^+(\m ,\b)}(\rho_b,\alpha)
\det'(1-b|(\m /\b)^+),
$$
and that for $\ga =a_\ga b_\ga$ split-regular we have
$$|W(\m_b,\b)|\prod_{\alpha \in \Phi_b^+(\m ,\b)}(\rho_b,\alpha)
\= |W(\g_\alpha ,\h)|\prod_{\alpha \in \Phi_\ga^+}(\rho_\ga
,\alpha),
$$
so that writing $L^M(\ga ,\tau) = L^M(b_\ga
  ,\tau)$ we get
$$
L^M(\ga ,\tau) \= \frac{\varpi_\ga( \det (1-\ga |
  (\k_M/\b)^+) \tr(\tau(b_\ga)))}
       {|W(\g_\ga ,\h)|\prod_{\alpha \in \Phi_\ga^+}(\rho_\ga
                        ,\alpha) \det'(1-\ga | (\m /\b)^+)}.
$$

From the above it follows that $\vol(\Ga_\ga \bs G_\ga) \CO_\ga
(f_t^0)$ equals:
$$
\frac{\chi_{_1}(X_\ga)l_{\ga_0}}
     {\det(1-\ga | \n)}\
\frac{e^{-l_\ga^2/4t}}
     {\sqrt{4\pi t}}\
a_\ga^{\rho_P} \sum_{c\ge 0,\ i\in I_c} (-1)^c\ e^{ts_i^c}\
\tr(b_\ga|\sigma_i^c)L^{M_2}(\ga ,\tau_i^c).
$$
Note that here we can replace $\sigma_i^c$ by
$\tilde{\sigma_i^c}$.

For a splitrank one Cartan $H=AB$ and a parabolic $P=MAN$ let
$\CE_P(\Ga)$ denote the set of $\Ga$-conjugacy classes
$[\ga]\subset\Ga$ such that $\ga$ is $G$-conjugate to an element
of $A^-B$, where $A^-$ is the negative Weyl chamber given by $P$.
We have proven:

\begin{theorem} \label{higher_heat_trace} Let $X_\Ga$ be a compact locally
Hermitian space with fundamental group $\Ga$ and such that the
universal covering is globally symmetric without compact factors.
Assume $\Ga$ is neat and write $\triangle_{p,q,\omega}$ for the
Hodge Laplacian on $(p,q)$-forms with values in a the flat
Hermitian bundle $E_\omega$, then it theta series defined by
 $$
 \Theta (t) \= \sum_{q=0}^{\dim_\C X_\Ga}
q(-1)^{q+1} \tr\ e^{-t\triangle_{0,q,\omega}}
 $$
equals
 $$
\sum_{P/conj.} \sum_{[\ga] \in \CE_P(\Ga)} \chi_{_1}(X_\ga)
\frac{l_{\ga_0}a_\ga^{\rho_P}} {\det(1-\ga | \n)}\
\frac{e^{-l_\ga^2/4t}} {\sqrt{4\pi t}}\
 $$ $$
\hspace{110pt}\times\sum_{c\ge 0,\ i\in I_c} (-1)^c\ e^{ts_i^c}\
\tr(b_\ga |\sigma_i^c)L^{M_2}(\ga ,\tau_i^c).
 $$ $$
        + f^0_t(e)\ \dim\omega\ vol(X_\Ga).\hspace{100pt}
 $$
\end{theorem}

The reader should keep in mind that by its definition we have for
the term of the identity:
$$
f^0_t(e)\ \dim\omega\ vol(X_\Ga) \= \sum_{q=0}^{\dim_\C
X}q(-1)^{q+1} \tr_\Ga(e^{-t\Delta_{0,q,\omega}}),
$$
where $\tr_\Ga$ is the $\Ga$-trace. Further note that by the
Plancherel theorem the Novikov-Shubin invariants of all operators
$\Delta_{0,q}$ are positive.

\subsection{The holomorphic torsion zeta function}

Now let $X$ be Hermitian again and let $\tau=\tau_1\otimes\tau_2$ be an irreducible representation of $K_{M}\rightarrow K_{M_1}\times K_{M_2}$.
Further assume that $\tau_1$ lies in the image of the restriction map $res: {\rm Rep}(M_1)\to{\rm Rep}(K_{M_1})$.

\begin{theorem}\label{hol_zeta}
Let $\Ga$ be neat and $(\omega ,V_\omega)$ a finite dimensional
unitary representation of $\Ga$. Choose a $\theta$-stable Cartan
$H$ of splitrank one. For $\Re(s)>>0$ define the zeta function
$Z_{P,\tau,\omega}^0(s)$ to be
$$
\exp \left( -\sum_{[\ga] \in \CE_H(\ga)} \frac{\chi_{_1}(X_\ga)
\tr(\omega(\ga)) \tr\tau_1(b_\ga)L^{M_2}(\ga ,\tau_2)}
     {\det (1-\ga|\n)}\
\frac{e^{-sl_\ga}}{\mu_\ga}\right).
$$
Then $Z_{P,\tau,\omega}^0$ has a meromorphic continuation to the
entire plane. The vanishing order of $Z_{P,\tau,\omega}^0(s)$ at
$s=\la (H_1)$, $\la \in \a^*$ is $(-1)^{\dim \n}$ times
$$
\sum_{\pi \in \hat{G}}N_{\Ga ,\omega}(\pi)
    \sum_{p,q,r}(-1)^{p+q+r} \dim \Big( H^q({\n},\pi_K)^\la
    \otimes \wedge^p\p_{M_1}\otimes\wedge^r\p_{M_2,-}
    \otimes V_{\breve{\tau}}\Big)^{K_M}.
$$
Note that
in the special case $M_2={1}$ this function equals the Selberg
zeta function.
\end{theorem}

\prf By Lemma 2.2.1 of \cite{lefschetz} the group $M_2$ is orientation preserving. 
Since $\tau_1$ lies in the image of the restriction,  the Euler-Poincar\'e function
$f_{\tau_1}^{M_1}$ for the representation $\tau_1$ exists. Further
for $M_2$ the function $g_{\tau_2}^{M_2}$ of Theorem 2.2.2 in \cite{lefschetz}
exists. We set $h_{\tau}(m_1,m_2) :=
f_{\tau_1}^{M_1}(m_1)g_{\tau_2}^{M_2}(m_2)$ and this function factors
over $M$. Then for any function $\eta$ on $M$ which is a product
$\eta =\eta_1\otimes\eta_2$ on $M_1\times M_2$ we have for the
orbital integrals:
$$
\CO_m^M(\eta) \= \CO_{m_1}^{M_1}(\eta_1)\CO_{m_2}^{{M_2}}(\eta_2).
$$
With this in mind it is straightforward to see that the proof of
Theorem \ref{hol_zeta} proceeds as the proof of Theorem
\ref{genSelberg} with the Euler-Poincar\'e function $f_\tau$
replaced by the function $h_{\tau}$.
\qed

Extend the definition of $Z_{P,\tau,\omega}^0(s)$ to arbitrary
virtual representations in the following way. Consider a finite
dimensional virtual representation $\xi = \oplus_i a_i \tau_i$
with $a_i\in \Z$ and $\tau_i\in \hat{K_M}$. Then let
$Z_{P,\xi,\omega}^0(s) = \prod_i Z_{P,\tau_i,\omega}^0(s)^{a_i}.$

\begin{theorem} \label{detformel}
Assume $\Ga$ is neat, then for $\la >>0$ we have the identity
$$
\prod_{q=0}^{\dim_\C X} \left(\frac{\det
(\Delta_{0,q,\omega}+\la)}{{\det}^{(2)}
(\Delta_{0,q,\omega}+\la)}\right)^{q(-1)^{q+1}}
            \=
$$ $$
\prod_{P/{\rm conj.}} \prod_{{c\ge 0}\ {i\in I_c}}
Z_{P,\sigma_i^c\otimes\tau_i^c ,\omega}(|\rho_P|+\sqrt{\la +
s_i^c})^{(-1)^c}
$$
\end{theorem}

\prf Consider Theorem \ref{higher_heat_trace}. For any
semipositive elliptic differential operator $D_\Ga$ the heat
trace $\tr e^{-tD_\Ga}$ has the same asymptotic as $t\to 0$ as
the $L^2$-heat trace $\tr_\Ga e^{-tD}$. Thus it follows that the
function
$$
h(t):= \sum_{q=0}^{\dim_\C X_\Ga} q(-1)^{q+1} (\tr
e^{-t\Delta_{0,q,\omega}}-\tr_\Ga e^{-t\Delta_{0,q,\omega}})
$$
is rapidly decreasing at $t=0$. Therefore, for $\la >0$ the Mellin
transform of $h(t)e^{-t\la}$ converges for any value of $s$ and
gives an entire function. Let
$$
\zeta_\la (s) := \frac 1{\Ga (s)} \int_0^\infty t^{s-1} h(t)
e^{-t\la} dt.
$$
We get that
$$
\exp(-\zeta_\la'(0)) \= \prod_{q=0}^{\dim_\C X} \left(\frac{\det
(\Delta_{0,q,\omega}+\la)}{{\det}^{(2)}
(\Delta_{0,q,\omega}+\la)}\right)^{q(-1)^{q+1}}.
$$
On the other hand, Theorem \ref{higher_heat_trace} gives a second
expression for $\zeta_\la(s)$. In this second expression we are
allowed to interchange integration and summation for $\la >>0$
since we already know the convergence of the Euler products
giving the right hand side of our claim.
\qed

Let $n_0$ be the order at $\la =0$ of the left hand side of the
last proposition. Then
$$
n_0 \= \sum_{q=0}^{\dim_\C(X)} q(-1)^q(h_{0,q,\omega}
-h_{0,q,\omega}^{(2)}),
$$
where $h_{0,q,\omega}$ is the $(0,q)$-th Hodge number of $X_\Ga$
with respect to $\omega$ and $h_{0,q,\omega}^{(2)}$ is the
$L^2$-analogue. Conjecturally we have $h_{0,q,\omega}^{(2)} =
h_{0,q,\omega}$, so $n_0=0$. For a splitrank one Cartan $H$, for
$c\ge 0$ and $i\in I_c$ we let
$$
n_{P,c,i,\omega} \ := \ {\rm
ord}_{s=|\rho_P|+\sqrt{s_i^c}}Z_{P,\sigma_i^c\otimes\tau_i^c
,\omega}(s)
$$
so $n_{P,c,i,\omega}$ equals $(-1)^{\dim \n}$ times
$$
\sum_{\pi \in \hat{G}}N_{\Ga ,\omega}(\pi)
    \sum_{p,q,r}(-1)^{p+q+r} \dim \Big( H^q({\n},\pi_K)^\la\otimes
\wedge^p\p_{M_1}\otimes\wedge^r\p_{M_2,-} \otimes
V_{\breve{\sigma}\otimes\breve{\tau}}\Big)^{K_M},
$$
for $\la(H)=|\rho_P|+\sqrt{s_i^c}$.
 We then consider
$$
c(X_\Ga ,\omega)\= \prod_{P}\prod_{c\ge 0,\ i\in I_c} \left(
2\sqrt{s_i^c}\right)^{(-1)^cn_{P,c,i,\omega}}.
$$

We assemble the results of this section to

\begin{theorem}
Let
$$
Z_\omega(s) \= \prod_{P/{\rm conj.}} \prod_{{c\ge 0}\ {i\in
I_c}}Z_{P,\sigma_i^c\otimes\tau_i^c,\omega}\left(
s+|\rho_P|+\sqrt{s_i^c}\right),
$$
then $Z_\omega$ extends to a meromorphic function on the plane.
Let $n_0$ be the order of $Z_\omega$ at zero then
$$
n_0 \= \sum_{q=0}^{\dim_\C(X)} q(-1)^q \left( h_{0,q}(X_\Ga)
-h_{0,q}^{(2)}(X_\Ga)\right),
$$
where $h_{p,q}(X_\Ga)$ is the $(p,q)$-th Hodge number of $X_\Ga$
and $h_{p,q}^{(2)}(X_\Ga)$ is the $(p,q)$-th $L^2$-Hodge number of
$X_\Ga$. Let $R_\omega(s) = Z_\omega(s)s^{-n_0}/c(X_\Ga ,\omega)$
then
$$
R_\omega(0) \= \frac{T_{hol}(X_\Ga
,\omega)}{T_{hol}^{(2)}(X_\Ga)^{\dim\omega}}.
$$
\end{theorem}

\prf This follows from Theorem \ref{detformel}. \qed

{\small Mathematisches Institut\\
Auf der Morgenstelle 10\\
72076 T\"ubingen\\
Germany\\
\tt deitmar@uni-tuebingen.de}

\today


\begin{thebibliography}{XXX}



\bibitem{AtSch}
\bf Atiyah, M.; Schmid, W.:
\it A Geometric Construction of the Discrete Series for Semisimple Lie Groups.
\rm Invent. math. 42, 1-62 (1977).


\bibitem{BM}
\bf Barbasch, D.; Moscovici, H.:
\it $L^2$-Index and the Selberg Trace Formula.
\rm J. Func. An. 53, 151-201 (1983).

\bibitem{BGV}
\bf Berline, N.; Getzler, E.; Vergne, M.:
\it Heat Kernels and Dirac Operators.
\rm Grundlehren 298. Springer 1992.


\bibitem{Bor}
\bf Borel, A.:
\it Introduction aux groupes arithm\'etiques.
\rm Hermann, Paris 1969.


\bibitem{BorWall}
 \bf Borel, A.; Wallach, N.:
 \it Continuous Cohomology, Discrete Groups, and Representations of Reductive Groups.
 \rm Ann. Math. Stud. 94, Princeton 1980.




\bibitem{BuOl}
\bf Bunke, U.; Olbrich, M.:
\it $\Gamma$-Cohomology and the Selberg Zeta Function
\rm J. reine u. angew. Math. 467, 199-219 (1995).





\bibitem{CloDel}
\bf Clozel, L.; Delorme, P.:
\it Le th\'eor\`eme de Paley-Wiener invariant pour les groupes de Lie
reductifs II.
\rm Ann. sci. \'Ec. Norm. Sup. (4) 23, 193-228 (1990).

\bibitem{D-Hitors}
\bf Deitmar, A.:
\it Higher torsion zeta functions.
\rm Adv. Math. 110, 109-128 (1995).





\bibitem{lefschetz}
\bf Deitmar, A.:
\it A higher rank Lefschetz formula. 
\rm J. Fixed Point Theory Appl. 2, 1-40 (2007).

%
%
\bibitem{Den}
\bf Deninger, C.:
\it A dynamical systems analogue of Lichtenbaum's conjectures on special values of Hasse-Weil zeta functions.
\rm http://arxiv.org/abs/math/0605724


\bibitem{DKV}
\bf Duistermaat, J.J.; Kolk, J.A.C.; Varadarajan, V.S.:
\it Spectra of locally symmetric manifolds of negative curvature.
\rm Invent. math. 52 (1979) 27-93.




\bibitem{Fried86}
\bf Fried, D.:
\it Analytic torsion and closed geodesics on hyperbolic manifolds.
\rm Invent. math. 84, 523-540 (1986).

\bibitem{Fried88}
\bf Fried, D.:
\it Torsion and closed geodesics on complex hyperbolic manifolds.
\rm Invent. Math. 91, 31-51 (1988).





\bibitem{GHJ}
\bf Goodman, F.M.; de la Harpe, P.; Jones, V.F.R.:
\it Coxeter Graphs and Towers of Algebras.
\rm Springer 1989.

\bibitem{GrSh}
\bf Gromov, M.; Shubin, M.A.:
\it Von Neumann Spectra near Zero.
\rm GAFA 1, 375-404 (1991).

\bibitem{HC-HA1}
\bf Harish-Chandra:
\it Harmonic analysis on real reductive groups I. The theory of
the constant term.
\rm J. Func. Anal. 19 (1975) 104-204.







\bibitem{HeSch}
\bf Hecht, H.; Schmid, W.:
\it Characters, asymptotics and $\n$-homology of Harish-Chandra modules.
\rm Acta Math. 151, 49-151 (1983).






\bibitem{Ju}
\bf Juhl, A.:
\it Zeta-Funktionen, Index-Theorie und hyperbolische Dynamik.\\
\rm Habilitationsschrift. Humboldt-Universit\"{a}t zu Berlin 1993.

\bibitem{Juhl}
 \bf Juhl, A.:
 \it Cohomological theory of dynamical zeta functions.
 \rm Progress in Mathematics, 194. Birkh\"auser Verlag, Basel,
2001.



\bibitem{Knapp}
\bf Knapp, A.:
\it Representation Theory of Semisimple Lie Groups.
\rm Princeton University Press 1986.



\bibitem{Lab}
\bf Labesse, J.P.:
\it Pseudo-coefficients tr\`es cuspidaux et K-th\'eorie.
\rm Math. Ann. 291, 607-616 (1991).


\bibitem{Licht}
\bf Lichtenbaum, S.:
\it The Weil-etale topology on schemes over finite fields.
\rm Compos. Math. 141, 689–702,  (2005).

\bibitem{L}
\bf Lott, J.:
\it Heat kernels on covering spaces and topological invariants.
\rm J.Diff. Geom. 35, 471-510 (1992)

\bibitem{LL}
\bf Lott, J.; L\"uck, W.:
\it $L^2$-Topologicla Invariants of 3-Manifolds.
\rm Invent. Math. 120, 15-60 (1995).


\bibitem{MS-eta}
\bf Moscovici, H.; Stanton, R.:
\it Eta invariants of Dirac operators on locally symmetric manifolds.
\rm Invent. math. 95, 629-666 (1989).

\bibitem{MS-tors}
\bf Moscovici, H.; Stanton, R.:
\it R-torsion and zeta functions for locally symmetric manifolds.
\rm Invent. math. 105, 185-216 (1991).


\bibitem{RS-RT}
\bf Ray, D.; Singer, I.:
\it R-torsion and the Laplacian on Riemannian manifolds.
\rm Adv. Math. 7, 145-220 (1971).




\bibitem{soule}
 \bf Soul\'e, C.; Abramovich, D.; Burnol, J.-F.; Kramer, J.:
 \it Lectures on Arakelov Geometry.
 \rm Cambridge University Press 1992.












\end{thebibliography}
\end{document}